\documentclass[11pt,reqno]{amsart}

\usepackage[T1]{fontenc}
\usepackage[utf8]{inputenc}
\usepackage{lmodern}
\usepackage{amsmath,amssymb,amsthm,mathtools}
\usepackage{mathrsfs}
\usepackage{enumitem}
\usepackage{microtype}
\usepackage{xcolor}
\usepackage[colorlinks=true,linkcolor=blue!55!black,citecolor=blue!55!black,urlcolor=blue!55!black]{hyperref}

\setlist[enumerate]{leftmargin=2.2em,itemsep=0.25em,topsep=0.35em}
\setlist[itemize]{leftmargin=2.0em,itemsep=0.25em,topsep=0.35em}
\numberwithin{equation}{section}

\newtheorem{theorem}{Theorem}[section]
\newtheorem{proposition}[theorem]{Proposition}
\newtheorem{lemma}[theorem]{Lemma}
\newtheorem{corollary}[theorem]{Corollary}

\theoremstyle{definition}
\newtheorem{definition}[theorem]{Definition}
\newtheorem{example}[theorem]{Example}

\newtheorem{conjecture}[theorem]{Conjecture}
\newtheorem{fact}[theorem]{Fact}
\theoremstyle{remark}
\newtheorem{remark}[theorem]{Remark}

\title [Proximal Myhill-type theorems]
{\textbf{Proximal Myhill-type theorems for $m$-expansive amenable semigroup actions}}

\author{Boxu Huang}
\address{School of Mathematics, Sichuan University, Chengdu 610064, Sichuan, China}
\email{mathhbx@163.com}

\subjclass[2020]{Primary 37B05, 37B40;  Secondary 20M14, 43A07}
\keywords{Garden of Eden theorem, Myhill property, amenable semigroup, $m$-expansiveness, weak specification, proximal pair, topological entropy}

\begin{document}

\begin{abstract}
The Garden of Eden theorem relates surjectivity to a weak injectivity condition for cellular automata on amenable groups. Li extended the Myhill implication to expansive actions of countable amenable groups with weak specification. We consider two simultaneous relaxations: the group is replaced by a cancellative amenable semigroup, and expansiveness by m-expansiveness. The loss of ordinary expansiveness removes homoclinic rigidity, so we introduce a proximal analogue of pre-injectivity. Our main result shows that, for a surjective action of an infinite countable cancellative left-amenable semigroup on a compact metric space, m-expansiveness plus weak specification forces every continuous equivariant map injective on each proximal set to be surjective. The proof uses entropy detection at a fixed scale, a strict entropy gap from weak specification, and finite multiplicity to produce proximal pairs from uniformly close points. We also give a one-dimensional variant under shadowing and examples clarifying the assumptions.
\end{abstract}

\maketitle

\section{Introduction}
The conceptual roots of the Garden of Eden theorem can be traced back to Richard Dedekind's classical characterization of finite sets in the 19th century: a set is finite if and only if every injective map from the set to itself is necessarily surjective. This elementary observation links finiteness with the functional property that injectivity implies surjectivity, planting a original of thought for later generations of dynamicists: if the state space of a system is an ``infinite product'' of finite sets, could there be some ``finite shadow'' that sustains an analogous injectivity-surjectivity relationship?

With the rise of symbolic dynamics, mathematicians began to study spaces of the form \(A^G\), where \(A\) is a finite alphabet and \(G\) is an infinite discrete group. Elements of \(A^G\) are called configurations, and, equipped with the product topology, they form a compact metric space. The group \(G\) acts on this space by translations: \((g\cdot x)(h)=x(g^{-1}h)\). A basic observation is that, although \(A^G\) is infinite, its configurations still retain a trace of 'finiteness'. Concretely, for any two distinct configurations \(x,y\in A^G\), the finiteness of the alphabet \(A\) guarantees the existence of some \(g\in G\) such that \(x\) and \(y\) differ at \(g\); in other words, the translation action distinguishes points. This reveals that, despite the vastness of the underlying space, the dynamics are rigidly constrained by the finite alphabet.

Building on these observations, W. Gottschalk formulated an interesting conjecture in the mid-20th century \cite{G}:
\begin{conjecture}[Gottschalk's Conjecture]
	For any finite set $A$ and any group $G$, if $T : A^G\to A^G$ is a continuous $G$-equivariant map(i.e., a cellular automaton), then injectivity of $T$ implies surjectivity. Here \(G\)-equivariance means that \(T(g\cdot x)=g\cdot T(x)\) for all \(g\in G\) and \(x\in A^G\).
\end{conjecture}
Groups satisfying this property are called \emph{surjunctive}. Gottschalk's conjecture weaves together Dedekind's finiteness insight, the rigidity phenomena of symbolic dynamics and the structure of abstract groups, making it one of the interesting problems at the interface of dynamical systems and group theory.

The Garden of Eden theorem is one of the basic rigidity principles in the theory of cellular automata.  For a finite alphabet $A$ and a countable group $G$, a cellular automaton is a continuous $G$-equivariant self-map of the full shift $A^G$.  The classical theorem of Moore and Myhill asserts, for $G=\mathbb Z^d$, that surjectivity is equivalent to pre-injectivity, where pre-injectivity means injectivity on configurations which differ only on a finite set \cite{Mo,My}.  The result was subsequently extended to cellular automata over amenable groups by Ceccherini-Silberstein, Mach\`i, and Scarabotti \cite{CMS}; the converse results of Bartholdi show that amenability is intrinsic to the classical Moore--Myhill phenomenon \cite{Ba1,Ba2}.  See \cite{CC} for a systematic account.

The dynamical content of the Myhill implication is particularly transparent.  A proper image should have strictly smaller orbit complexity, whereas pre-injectivity prevents the collapse of orbit patterns which differ only ``at finite distance''.  In symbolic systems these two facts can be made combinatorial.  For general compact dynamical systems the appropriate replacement for finite disagreement is the homoclinic, or asymptotic, relation.  In this direction, Li proved that every expansive action of a countable amenable group with the weak specification property has the Myhill property: every continuous equivariant map which is injective on each homoclinic class is surjective \cite{Li}.  The proof is driven by a strict entropy-gap theorem for proper invariant subsystems and by a mechanism which forces an entropy-losing factor map to identify two homoclinic points.

The purpose of the present paper is to examine how far that mechanism survives when two of its structural hypotheses are weakened simultaneously.  First, we replace the acting group by a cancellative amenable semigroup.  This is not a formal change.  In a semigroup there are no inverses with which to move orbit windows freely, left and right translations play different roles, and amenability is naturally one-sided.  The filling theorem of Ceccherini-Silberstein, Coornaert, and Krieger provides the appropriate substitute for the Ornstein--Weiss quasitiling machinery in the cancellative-semigroup setting \cite{CCK}.  Second, we replace expansiveness by $m$-expansiveness, introduced by Morales \cite{Mora}.  Here an infinite Bowen ball is allowed to contain finitely many, but at most $m$, points.  Thus $m$-expansiveness retains finite orbit multiplicity while losing the uniqueness built into classical expansiveness.

This loss of uniqueness is precisely where Li's homoclinic argument ceases to be robust \cite{Li}. For an expansive action, sufficiently strong uniform closeness of two complete orbits forces equality, and uniform closeness outside a finite set can be upgraded to asymptoticity.  Neither implication is available for a general $m$-expansive action.  A finite collection of distinct points may remain uniformly close for all time without becoming asymptotic.  The natural replacement is therefore proximality.  Recall that $x,y\in X$ are proximal when
\[
\inf_{s\in S} d(sx,sy)=0.
\]
Compactness and finite orbit multiplicity interact well with this weaker relation: if more than $m$ distinct points remain inside the same $m$-expansive scale outside a finite set, then at least two of them must be proximal.  This elementary observation is the bridge which allows the entropy argument to survive.

The second new issue is entropy at a fixed scale.  Li's proof uses the fact that expansiveness allows topological entropy to be computed from separated sets at one sufficiently small scale.  For $m$-expansive semigroup actions this is not immediate.  We prove a fixed-scale entropy theorem by combining the finite cardinality of infinite Bowen balls with the filling theorem for cancellative amenable semigroups.  The key point is that each local infinite Bowen ball has at most $m$ possible ``cores'', while the number of quasitile centers can be made arbitrarily small relative to the volume of a F{\o}lner set.  Consequently the local multiplicity contributes zero entropy.

Our main theorem is the following proximal version of the Myhill implication.

\begin{theorem}\label{main}
	Let $S$ be an infinite countable cancellative left-amenable semigroup acting continuously on a compact metric space $(X,d)$.  Assume that every action map $x\mapsto sx$ is surjective, that the action is $m$-expansive, and that it has the weak specification property.  Let $f\colon X\to X$ be continuous and $S$-equivariant.  If
	\[
	f|_{{\rm{Prox}}(x)} \text{ is injective for every }x\in X,
	\]
	then $f$ is surjective.
\end{theorem}

The formulation deserves two comments.  First, the proximal relation need not be transitive, so $\rm{Prox}(x)$ is a proximal set rather than an equivalence class.  Second, injectivity on every proximal set is stronger than the usual pairwise requirement that two proximal points with the same image coincide.  The stronger formulation is the one naturally compatible with the common-proximal-center phenomenon produced by weak specification; see Proposition~\ref{prop:proximal-dense-Gdelta} and Corollary~\ref{0202} below.

The proof of Theorem~\ref{main} is organized around two structural results.  The first is an entropy gap.

\begin{theorem}[Entropy gap]\label{thm:intro-gap}
	Under the hypotheses on $(X,S)$ in Theorem~\ref{main}, every nonempty proper closed $S$-invariant subset $Y\subsetneq X$ satisfies
	\[
	h_{\rm {top}}(Y)<h_{\rm {top}}(X).
	\]
\end{theorem}

The second says that an equivariant factor which loses entropy must collapse a proximal pair.

\begin{theorem}[Proximal collision]\label{thm:intro-collision}
	Let $(X,S)$ satisfy the hypotheses of Theorem~\ref{main}, let $f\colon X\to X$ be continuous and $S$-equivariant, and put $Y=f(X)$.  If
	\[
	h_{\rm {top}}(Y)<h_{\rm {top}}(X),
	\]
	then there exist distinct proximal points $x_1,x_2\in X$ with
	\[
	f(x_1)=f(x_2).
	\]
\end{theorem}

Theorem~\ref{main} follows immediately: if $f(X)$ were proper, Theorem~\ref{thm:intro-gap} would give a strict entropy loss, and Theorem~\ref{thm:intro-collision} would contradict proximal-set injectivity.

For a single map, shadowing restores enough rigidity to recover an asymptotic formulation.  We prove the following complementary result.

\begin{theorem}\label{thm:one-dimensional-main}
	Let $T\colon X\to X$ be a surjective continuous map of a compact metric space.  Assume that the $\mathbb N_0$-action generated by $T$ is $m$-expansive and has the weak specification property, and assume moreover that $T$ has the shadowing property.  If $f\colon X\to X$ is continuous, commutes with $T$, and is injective on every forward asymptotic class, then $f$ is surjective.
\end{theorem}

The shadowing hypothesis in Theorem~\ref{thm:one-dimensional-main} is essential for the argument: it is used to show that a sufficiently small local stable set contains at most $m$ forward asymptotic classes.  Without shadowing, $m$-expansiveness alone does not provide this zero-scale asymptotic rigidity.

The paper is organized as follows.  Section 1 collects the semigroup and entropy preliminaries.  In Section 2 we establish the fixed-scale entropy formula for $m$-expansive cancellative amenable semigroup actions.  Section 3 develops the proximal tools, including residuality of proximal sets for infinite cancellative monoids with weak specification.  Section 4 proves the entropy gap, the proximal collision theorem, and Theorem~\ref{main}.  Section 5 proves Theorem~\ref{thm:one-dimensional-main}.  Section 6 contains examples showing that several hypotheses cannot simply be removed.


\section{Preliminaries}

In this section, we introduce the notation and fundamental concepts and review some background  material on dynamical systems, amenable semigroups, topological entropy, weak specification property, $m$-expansive system and proximal pair.

\subsection{General notation.}

Throughout this paper, we denote by $|\cdot|$ the cardinality of finite sets and $S$ denote by a countably infinite semigroup.

Given a semigroup $S$, we use ${\rm{Fin}}(S)$ to denote the set of all nonempty finite subsets of $S$.

Given a metric space $(X, d)$ and a subset $Z\subseteq X$,  we define the {\it {diameter}} of $Z$ as ${\rm {diam}}(Z, d):=\sup \{d(x_1, x_2) \, :\, x_1, x_2\in Z\}$.
The {\it{distance}} between two subsets $A, B\subseteq X$ is denoted by
$$
{\rm {dist}}_{d}(A, B):=\inf_{x\in A, y\in B}\{d(x, y)\}.
$$
For $x\in X$ and $r>0$ we denote $B_d(x, r)$ (resp. $\overline{B}_d(x, r)$) the open ball (resp. the close ball) centered at $x$ with radius $r$ by
\begin{equation*}
	B_{d}(x, r):=\{y\in X \, :\, d(x, y)<r\}\quad \left(\text {resp.  }\   \overline{B}_{d}(x, r):=\{y\in X \, :\, d(x, y)\leq r\}\right).
\end{equation*}

\subsection{Cancellative amenable semigroups}

Recall that a {\it {semigroup}} is a set equipped with an associative binary operation and a {\it {monoid}} is a semigroup admitting an identity element. Let $S$ be a countable infinite semigroup. We use a multiplicative notation for the binary operation on $S$.

Given $s\in S$, we use $L_s$ and $R_s$ to denote the left and right multiplication by $s$, that is, the maps $L_s : S\rightarrow S$ and $R_s : S\rightarrow S$ defined
by
\begin{equation*}
	L_s(s')=s s'\quad {\text {and \  }}\  R_s(s')=s's\quad {\text {for all \  }}\   s'\in S.
\end{equation*}

\begin{definition}\label{cancel}
	An element $s\in S$ is called {\it {left-cancellative}}(resp. {\it {right-cancellative}}) if the map $L_s$ (resp. $R_s$)
	is injective. We say that an element $s\in S$ is {\it {cancellative}} if it is both left-cancellative and right-cancellative. The semigroup is called {\it {left-cancellative}}
	(resp. {\it {right-cancellative}}, resp. {\it {cancellative}}) if every element in $S$ is {\it {left-cancellative}}
	(resp. {\it {right-cancellative}}, resp. {\it {cancellative}}).
\end{definition}

We also record a simple finite-preimage estimate which will be used later.

\begin{fact}\label{finite set lemma}
	If $S$ is left cancellative and $K, W\subset S$ are finite, then there exists a finite set $E\subset S$ such that
	\[
	W(S\setminus E)\subset S\setminus K.
	\]
\end{fact}

\begin{proof}
	Take
	\[
	E=\bigcup_{w\in W}L_w^{-1}(K).
	\]
	Since each $L_w$ is injective, $|L_w^{-1}(K)|\le |K|$, hence $E$ is finite.  The desired inclusion is immediate from the definition.
\end{proof}

We next formalize the notion of approximate invariance of one finite subset of the semigroup $S$ with respect to another.

\begin{definition}
	Let $K$ and $A$ be nonempty finite subsets of a semigroup $S$. We define the $K$-{\it {interior}} ${\rm {Int}}_K(A)$ and the $K$-{\it {boundary}} $\partial_K(A)$ of
	$A$ by
	\begin{equation*}
		\begin{split}
			{\rm {Int}}_K(A) & :=\{s\in A : Ks\subset A\}, \\
			\partial_K(A)  & :=A\backslash  {\rm {Int}}_K(A).
		\end{split}
	\end{equation*}
\end{definition}

From the definitions, it is easy to see the following proposition (\cite[Proposition 2.1]{CCK}) holds.

\begin{proposition}\label{hh01}
	Let $S$ be a semigroup and let $K, A, B$ be finite subsets of $S$. Then
	\begin{itemize}
		\item[(1)] if every element in $K$ is left-cancellative, then $\partial_K(A)=\bigcup_{k\in K} L_k^{-1}(k A\backslash A)$;
		\item[(2)] $\partial_K(A\cup B)\subset\partial_K(A)\cup \partial_K(A)$;
	\end{itemize}
\end{proposition}

\begin{definition}
	Let $S$ be a left-cancellative semigroup. We then define the {\it {relative amenability constat}} $\alpha(A, K)$ of $A$ with respect to $K$ by
	\begin{equation*}
		\alpha(A, K):=\frac{| \partial_K(A)|}{|A|}.
	\end{equation*}
\end{definition}

In this paper, we need the following basic property of the relative amenability constat.

\begin{proposition}\label{racp}
	Let $0<\theta<1/2$, $S$ be a left cancellative semigroup and let $F, B, K\in {\rm {Fin}}(S)$. Suppose that
	$B\subset F$, $|F\backslash B|\leq \theta|F|$ and $\alpha(F, K)\leq\theta$. Then
	\begin{equation*}
		\alpha(B, K)\leq 4|K|\theta.
	\end{equation*}
\end{proposition}

\begin{proof}
	Denote by
	\begin{equation*}
		T=\bigcup_{s\in K} L_s^{-1}(F\backslash B).
	\end{equation*}
	Suppose that $g\in {\rm {Int}}_K(F)\cap B$, then $g\in B$ and $K g\subseteq F$. If $K g\cap (F\backslash B)\not=\emptyset$, then
	\begin{equation*}
		g\in \bigcup_{s\in K} L_s^{-1}(F\backslash B)=T.
	\end{equation*}
	Thus we deduce that
	\begin{equation}\label{inteq01}
		{\rm {Int}}_K(B)\supseteq ({\rm {Int}}_K(F)\cap B)\backslash T.
	\end{equation}
	Meanwhile, since $|F\backslash B|\leq \theta|F|$ and $\alpha(F, K)\leq\theta$, it is easy to see that
	\begin{equation}\label{inteq02}
		|{\rm {Int}}_K(F)\cap B|\geq |{\rm {Int}}_K(F)|-|F\backslash B|\geq (1-2\theta)|F|.
	\end{equation}
	Furthermore, owing to $S$ being left cancellative semigroup, we get
	\begin{equation}\label{inteq03}
		|T|\leq \sum_{s\in K}|L_s^{-1}(F\backslash B)|\leq \theta|K| |F|.
	\end{equation}
	Note that $B\subset F$ and $|F\backslash B|\leq \theta|F|$ implies that $\frac{|F|}{|B|}\geq\frac{1}{1-\theta}$.
	Combining with (\ref{inteq01}), (\ref{inteq02}) and (\ref{inteq03}), we deduce that
	\begin{equation*}
		\frac{|{\rm {Int}}_K(B)|}{|B|}\geq (1-(2+|K|)\theta)\frac{|F|}{|B|}\geq\frac{1-(2+|K|)\theta}{1-\theta}.
	\end{equation*}
	Hence, noting $\theta<1/2$, we get
	\begin{equation*}
		\alpha(B, K)=\frac{|B|-|{\rm {Int}}_K(B)|}{|B|}\leq \frac{(1+|K|)\theta}{1-\theta}\leq  2(1+|K|)\theta\leq 4|K|\theta.
	\end{equation*}
	Therefore, the proof of the proposition is completed.
\end{proof}

\begin{definition}
	Let $X$ be a set and $\varepsilon>0$. We say that a collection $\{A_i\}_{i\in I}$ of finite subsets of $X$ is $\varepsilon$-{\it {disjoint}} if there exist pairwise disjoint sets $\widehat{A_i}\subseteq A_i$ such that $|\widehat{A_i}|\geq (1-\varepsilon)|A_i|$ for all $i\in I$.
\end{definition}

\begin{definition}
	Let $S$ be a semigroup. Let $K$ and $A$ be finite subsets of $S$. For $\varepsilon>0$, a finite subset $P\subset S$ is called an $(\varepsilon, K)$-{\it {filling pattern}} for $A$ if the following conditions are satisfied:
	\begin{itemize}
		\item[(1)] $P\subset{\rm {Int}}_K(A)$;
		\item[(2)] the collection $\{Ks\}_{s\in P}$ is $\varepsilon$-disjoint.
	\end{itemize}
\end{definition}

In this paper, we need the following key tool, filling theorem \cite[Theorem 3.8]{CCK}, which is a semigroup analogue of the famous
the Ornstein-Weiss quasitiling theorem for groups.

\begin{theorem}[Filling theorem]\label{Filling theorem}
	Let $S$ be a cancellative semigroup and let $0<\varepsilon\leq 1/2$. There exists an integer $N\geq 1$, depending only on $\varepsilon$, such that the following holds.
	
	If $\{K_i\}_{1\leq i\leq N}$ is a finite sequence of nonempty finite subsets of $S$ with
	\begin{equation}\label{filling01}
		\alpha(K_j, K_i)\leq \frac{\varepsilon^{2N}}{|K_i|} \quad {\text {for all}}\quad 1\leq i<j\leq N
	\end{equation}
	and $F$ is a nonempty finite subset of $S$ such that
	\begin{equation}\label{filling02}
		\alpha(F, K_i)\leq \varepsilon^{2N}  \quad {\text {for all}}\quad 1\leq i\leq N,
	\end{equation}
	then there exists a finite sequence $\{P_i\}_{1\leq i\leq N}$ of finite subsets of $S$ such that
	\begin{itemize}
		\item[(1)] each $P_i, 1\leq i\leq N$, is an $(\varepsilon, K_i)$-filling pattern of $F$;
		\item[(2)] the subsets $K_iP_i\subset F$, $1\leq i\leq N$, are pairwise disjoint;
		\item[(3)] $F^*:=F\backslash \bigcup_{1\leq i\leq N} K_iP_i$ satisfies $|F^*|\leq \varepsilon |F|$.
	\end{itemize}
\end{theorem}

\subsection{Amenable semigroups.}

This section delves into the fundamental definitions and characterizations of amenable semigroups and in the process introduce some of the basic notation and terminology appearing throughout this paper.

A left-invariant mean on $S$ is a positive linear functional $\mathfrak m\colon \ell^\infty(S)\to\mathbb R$ satisfying $\mathfrak m(1)=1$ and
\[
\mathfrak m(\varphi\circ L_s)=\mathfrak m(\varphi)
\qquad(s\in S,\ \varphi\in\ell^\infty(S)).
\]
The semigroup is \emph{left amenable} if it admits a left-invariant mean.

\begin{definition}\label{def-Folner}
	Let $S$ be a semigroup. A \emph{left-invariant mean} on $S$ is a positive linear functional $\mathfrak m\colon \ell^\infty(S)\to\mathbb R$ satisfying $\mathfrak m(1)=1$ and
	\[
	\mathfrak m(\varphi\circ L_s)=\mathfrak m(\varphi)
	\qquad(s\in S,\ \varphi\in\ell^\infty(S)).
	\]
	The semigroup is \emph{left amenable} if it admits a left-invariant mean.
\end{definition}

The notion of an amenable group was introduced in 1929 by von Newmann \cite{Ne} for the study of the Banach-Tarski paradox. Subsequently, the conception of  amenable semigroup was introduced in the 1940s by Day; see \cite{D1,D2,D3}. In particular, Day \cite{D1} proved that every commutative semigroup is amenable which generalized a results previously obtained by von Newmann for groups. In fact, compared with amenable groups, amenable semigroups have many new properties that amenable groups do not have. For example, left-amenable and right-amenable are equivalent for groups, every finite group is amenable, and every subgroup of amenable group is amenable.
However, there exists finite semigroup that is left-amenable but not right-amenable and finite amenable semigroup containing sub-semigroup that is neither
left amenable nor right amenable.

\begin{definition}
	A semigroup $S$ satisfies {\it {the F{\o}ner condition}} (FC) if for every finite subset $K\subset S$ and every $\varepsilon>0$, there exists a nonempty
	finite subset $F\subseteq S$ such that
	\begin{equation}\label{FC-condition}
		|k F\backslash F|\leq \varepsilon |F| \quad {\text {for all \  }}\ k\in K.
	\end{equation}
\end{definition}

In his thesis, Frey \cite{Fr} showed that every left-amenable semigroup $S$ satisfies the F{\o}ner condition (FC); see \cite[Theorem 3.5]{Na} for a simpler proof.
However, the F{\o}ner condition (FC) is a necessary but not sufficient condition for left-amenability of semigroups. Examples of semigroups that are not left-amenable
but satisfy the F{\o}ner condition (FC) are provided by finite semigroup $S$ that are not left-amenable. In fact, it is clear that taking $F=S$ shows that any finite
semigroup $S$ trivially satisfies the F{\o}ner condition (FC).

\begin{definition}
	Let $S$ be a countable semigroup. A sequence $\{F_n\}_{n\in\mathbb{N}}$ of nonempty finite subsets of $S$ is said to be {\it {left F{\o}ner sequence}} if for
	every $s\in S$ we have that
	\begin{equation}\label{FS}
		\lim_{n\rightarrow\infty} \frac{|s F_n\backslash F_n|}{|F_n|}=0.
	\end{equation}
\end{definition}

\begin{remark}\label{amenable_remark}
	Let $S$ be a countably left amenable semigroup. A sequence of finite subsets $\{F_n\}$ of $S$
	is a left F{\o}ner sequence if and only if it satisfies that for any $s\in S$,
	\begin{equation*}
		\lim_{n\rightarrow\infty} \frac{|s F_n\cap F_n|}{|F_n|}=1.
	\end{equation*}
\end{remark}

\begin{proposition}\label{FC-FS}
	A countable semigroup $S$ satisfies the F{\o}ner condition (FC) if and only if it admits a left F{\o}ner sequence.
\end{proposition}

\begin{proof}
	Indeed, if $S$ satisfies condition (FC), we can construct a sequence $\{F_n\}_{n\in\mathbb{N}}$ satisfying (\ref{FS}) in the following way. We first choose
	an increasing sequence $S_1\subseteq S_2\cdots $ of finite subsets of $S$ whose union is $S$, that is,
	\begin{equation*}
		S=\bigcup_{n=1}^{\infty} S_n, \quad S_n\in{\rm {Fin}}(S) \quad {\text {and}} \quad S_n\subset S_{n+1}.
	\end{equation*}
	Then, for every finite subset $S_n\subset S$ and every $\varepsilon>0$, we take as $F_n$ one of the nonempty finite subsets of $S$ satisfying (\ref{FC-condition}).
	Hence $\{F_n\}_{n\in\mathbb{N}}$ is a left F{\o}ner sequence.
	
	Conversely, suppose that the $\{F_n\}_{n\in\mathbb{N}}$ is a left F{\o}ner sequence. Let $K\subset S$ be a finite subset and $\varepsilon>0$. Then, for every
	$k\in K$, there is $i_k\in\mathbb{N}$ such that $|k F_i\backslash F_i|\leq \varepsilon |F_i|$ for all $i\geq i_k$. To get a nonempty finite subset $F\subset S$
	satisfying (\ref{FC-condition}), it suffices to take $F:=F_i$, where $i\in\mathbb{N}$ satisfies that $i\geq i_k$ for every $k\in K$.
\end{proof}

For the left-cancellative semigroup, combining with above arguments, we know that the following characterizations of left-amenability holds (also see \cite[Theorem 3.1]{CC1}).

\begin{theorem}[F{\o}ner-Frey-Namioka]
	Let $S$ be a countable left-cancellative semigroup. Then the following statements are equivalent:
	\begin{itemize}
		\item[(1)] $S$ is left-amenable;
		\item[(2)] $S$ satisfies the F{\o}ner condition (FC);
		\item[(3)] $S$ admits a left F{\o}ner sequence.
	\end{itemize}
\end{theorem}

We shall use the following well known results for amenable groups about F{\o}ner sequence.

\begin{proposition}\label{infinite}
	Let $S$ be a countable left-cancellative and left amenable semigroup and let $\{F_n\}_{n\in \mathbb{N}}$ be a left F{\o}ner sequence for $S$. Let $K$ be a nonempty subset of $S$, $\eta>0$ and $N\in\mathbb{N}$. Then
	the following hold:
	\begin{itemize}
		\item[(1)] one has
		\begin{equation*}
			\lim_{n\rightarrow\infty} |F_n|=\infty;
		\end{equation*}
		\item[(2)] one has
		\begin{equation*}
			\lim_{n\rightarrow\infty} \frac{|K F_n\backslash F_n|}{|F_n|}=0\quad {\text {and}}\quad \lim_{n\rightarrow\infty} \alpha(F_n, K)=0.
		\end{equation*}
		\item[(3)] there exist nonempty finite subsets $K_1, \cdots, K_N$ of $S$ satisfying: $K\subset K_i$ for all $1\leq i\leq N$ and
		\begin{equation}\label{amenable01}
			\alpha(K_j, K_i)\leq \frac{\eta^{2N}}{|K_i|}\quad {\text {for all}}\quad 1\leq i<j\leq N.
		\end{equation}
	\end{itemize}
\end{proposition}

\begin{proof}
	Let $k\in\mathbb{N}$. Since $S$ is infinite, we can choose a finite subset $E\subset S$ with cardinality $|E|=k^2$. It follows from (\ref{FS}) that, for each $s\in E$,
	there is $N(s)\in\mathbb{N}$ such that $|s F_n\backslash F_n|<|F_n|$ for all $n\geq N(s)$. The left-cancellativity of $S$ means that $|s F_n|=|F_n|$ and
	hence $s F_n$ meets $F_n$ for all  $n\geq N(s)$. Again, based on the left-cancellativity of the semigroup $S$, we can deduce that, for $n\geq N(s)$,
	\begin{equation*}
		s\in\bigcup_{g\in F_n} L_g^{-1}(F_n)\quad {\text {and}}\quad \left|\bigcup_{g\in F_n} L_g^{-1}(F_n)\right|\leq |F_n|^2.
	\end{equation*}
	If we take $N:=\max_{s\in E} N(s)$, we then obtain that $E\subseteq \cup_{g\in F_n} L_g^{-1}(F_n)$ for all $n\geq N$. Thus we get that
	\begin{equation*}
		k^2=|E|\leq |F_n|^2
	\end{equation*}
	and hence $k\leq |F_n|$ for all $n\geq N$. This means that (1) holds.
	
	To prove (2), we first observe that
	\begin{equation*}
		\begin{split}
			\frac{|K F_n\backslash F_n|}{|F_n|} & =\frac{\left|\bigcup_{s\in K} (s F_n\backslash F_n)\right|}{|F_n|} \\
			& \leq \sum_{s\in K} \frac{|s F_n\backslash F_n|}{|F_n|}.
		\end{split}
	\end{equation*}
	This implies the first equality of (2) since the right side of above inequality tends to $0$ as $n\rightarrow\infty$ by (\ref{FS}). Meanwhile, by Proposition \ref{hh01} and the definition of left-cancellative semigroup, it follows that
	\begin{equation*}
		\alpha(F_n, K)=\frac{|\partial_K(F_n)|}{|F_n|}\leq \sum_{k\in K}\frac{k F_n\backslash F_n}{|F_n|}.
	\end{equation*}
	This shows that the second equality of (2) holds.
	
	To prove (3), we construct recursively finite subsets $K_1, K_2,\cdots, K_N$ of $S$ which have the required properties. Set $K_1:=K$. Suppose that the $K_1, \cdots, K_k$
	has been constructed recursively with $1\leq k\leq N-1$ which satisfy the properties (\ref{amenable01}). By (1) and (2) we know that
	\begin{equation*}
		\lim_{n\rightarrow\infty}\frac{|K|}{|F_n|}=0 \quad {\text {and}}\quad \lim_{n\rightarrow\infty}\alpha(F_n, K_i)=0
		\quad {\text {for all}}\quad 1\leq i\leq k.
	\end{equation*}
	Thus there exists $F_p\in\{F_n\}_{n\in\mathbb{N}}$ satisfies that
	\begin{equation*}
		\frac{|K|}{|F_p|}<\frac{\eta^{2N}}{2\max_{1\leq j\leq k}|K_j|} \quad {\text {and}}\quad \alpha(F_p, K_i)<\frac{\eta^{2N}}{2\max_{1\leq j\leq k}|K_j|}
		\quad {\text {for all}}\quad 1\leq i\leq k.
	\end{equation*}
	Take $K_{k+1}:=F_p\cup K$. By Proposition \ref{hh01}, it follows that, for $1\leq i\leq k$,
	\begin{equation*}
		\begin{split}
			\alpha(K_{k+1}, K_i) & =\frac{|\partial_{K_i}(K_{k+1})|}{|K_{k+1}|} \leq \frac{|\partial_{K_i}(F_{p})|}{|F_p|}+\frac{|\partial_{K_i}(K)|}{|F_p|} \\
			& \leq \alpha(F_p, K_i)+\frac{|K|}{|F_p|}\\
			& <\frac{\eta^{2N}}{\max_{1\leq j\leq k}|K_j|}\leq\frac{\eta^{2N}}{|K_i|}.
		\end{split}
	\end{equation*}
	Hence the proof of the proposition is completed.
\end{proof}

\begin{theorem}\label{foner111}
	Let $S$ be a countable left amenable and right-cancellative semigroup. Then for any finite subsets $F$, $K$ of $S$ and $0<k<1$, there exists a finite subsets $B$ of
	$S$ such that $B\supseteq K$ and
	\begin{equation*}
		|s B\cap B|\geq k |B| \quad {\text {for each  }}\   s\in F.
	\end{equation*}
\end{theorem}

\begin{proof}
	Choose $n$ large enough with
	\begin{equation}\label{fon01}
		k\left(1+\frac{|K|}{n}\right)<1.
	\end{equation}
	
	In what follows we wish to construct a finite subset $C$ of $S$ satisfying
	\begin{enumerate}
		\item $|C|\geq n$, and
		\item $|s C\cap C|\geq k(1+|K|/n) |C|$ for each $s\in F$.
	\end{enumerate}
	This is achieved by means of the following construction.
	Let $F^*$ be a finite subset of $S$ such that $F^*\supseteq F$ and $|F^*|\geq n^2$ since $S$ is an infinite semigroup.
	Since $S$ is left amenable and Remark \ref{amenable_remark}, there exists a left-F{\o}ner sequence $\{F_m\}_{m\in\mathbb{N}}$ for $S$ which satisfies that
	\begin{equation*}
		\lim_{m\rightarrow\infty} \frac{|s F_m\cap F_m|}{|F_m|}=1 \quad \text{for all }\ s\in F^*.
	\end{equation*}
	It follows that there must be some $m\in\mathbb{N}$ with
	\begin{equation*}
		|s F_m\cap F_m|\geq k\left(1+\frac{|K|}{n}\right) |F_m|
	\end{equation*}
	for each $s\in F^*$ (hence in $F$). Take the finite set $C$ of $S$ to be such $F_m$ above.
	This finite set $C$, therefore, satisfies the property (2). To see (1), notice that $|s C\cap C|>0$ (that is
	$s C\cap C\not=\emptyset$) for each $s\in F^*$. This implies that
	\begin{equation*}
		s\in \bigcup_{c'\in C} R_{c'}^{-1}(C), \quad {\text {for each  }}\  s\in F^*.
	\end{equation*}
	Here  $R_{c'} : S\rightarrow S$ is the right-multiplication by $c'$, that is, the map defined by $R_{c'}(g)= g c'$ for all $g\in S$. Since the
	map $R_{c'}$ is injective ($S$ being right-cancellative), one has $|R_{c'}^{-1}(C)|\leq |C|$. Hence we get
	\begin{equation*}
		n^2\leq|F^*| \leq \left|\bigcup_{c'\in C} R_{c'}^{-1}(C)\right|\leq |C|^2,
	\end{equation*}
	showing that $|C|\geq n$ holds.
	
	Now it remains to produce a finite subset $B$ of $S$ satisfying the statements of the theorem.
	Let $B=C\cup K$, then obviously $B$ is a finite subset of $S$ and, for each $s\in F$,
	\begin{equation*}
		\begin{split}
			|s B\cap B|\geq |s C\cap C| & \geq k\left(1+\frac{|K|}{n}\right) |C| \\
			& \geq k\left(1+\frac{|K|}{|C|}\right)|C|\\
			& =k(|C|+|K|)\geq k |B|.
		\end{split}
	\end{equation*}
	Hence the theorem is proved.
\end{proof}

\begin{theorem}\label{foner_exist}
	Let $S$ be a countable left amenable and right-cancellative semigroup. Then there exists a left F{\o}ner sequence $\{F_n\}$ of $S$ satisfying
	$F_n\subseteq F_{n+1}$ and $\bigcup_{n=1}^{\infty} F_n=S$.
\end{theorem}

\begin{proof}
	Fix an increasing sequence $E_1\subseteq E_2\cdots $ of finite subsets of $S$ whose union is $S$. We will construct, recursively over $n$, nonempty finite subsets
	$F_n$ of $S$ with the following properties:
	\begin{enumerate}
		\item for every $n\in\mathbb{N}$, $E_n\subseteq F_n$
		\item for every $n\in\mathbb{N}$, one has
		\begin{equation*}
			|s F_n\cap F_n|\geq \left(1-\frac{1}{n}\right) |F_n| \quad {\text {for all  }}\  s\in E_n.
		\end{equation*}
	\end{enumerate}
	
	Let $0<k_1<1$. By Theorem \ref{foner111}, there exists
	a finite subset $F_1$ of $S$ satisfying $F_1\supseteq E_1$ and $|s F_1\cap F_1|\geq k_1 |F_1|$ for all $s\in E_1$. Now let $n\geq 1$ and suppose that
	we have constructed nonempty finite subsets $F_n$ of $S$ which satisfy the required properties. Set $F=E_{n+1}$, $K=F_n\cup E_{n+1}$
	and $k=\left(1-\frac{1}{n+1}\right)$ as in Theorem \ref{foner111}. Then by Theorem \ref{foner111}
	we can find a nonempty finite subset $F_{n+1}$ of $S$ such that
	$F_{n+1}\supseteq F_n\cup E_{n+1}$ and
	\begin{equation*}
		|s F_{n+1}\cap F_{n+1}|\geq \left(1-\frac{1}{n+1}\right) |F_{n+1}| \quad {\text {for all  }}\   s\in E_{n+1}.
	\end{equation*}
	Then properties (1) and (2) hold for $n+1$, which completes the recursive step. Let $g\in S$. Then $g\in E_n$ for large enough $n$. Thus we have
	\begin{equation*}
		\lim_{n\rightarrow\infty}\frac{|g F_n\cap F_n|}{|F_n|}=1.
	\end{equation*}
	By Remark \ref{amenable_remark}, we know that $\{F_n\}$ is a left-F{\o}ner sequence of $S$.
	Hence $\{F_n\}$ is a left-F{\o}ner sequence of $S$ satisfying
	$F_n\subseteq F_{n+1}$ and $\bigcup_{n=1}^{\infty} F_n=S$ as desired.
\end{proof}

In this paper, we need  an analogue of Fekete lemma for subadditive and right-subinvariant functions defined on the set of all finite subsets of a cancellative left amenable semigroup. A suggested reference for the following theorem is Ceccherini-Silberstein, Coornaert and Krieger's paper \cite[Theorem 1.1]{CCK}.
\begin{theorem}[\cite{CCK}]\label{fekete}
	Let $S$ be a countably cancellative left amenable semigroup and let $h : {\rm{Fin}}(S)\rightarrow \mathbb{R}$ be a real-valued map satisfying the
	following conditions:
	\begin{itemize}
		\item[(1)] $h$ is subadditive, i.e.,
		$$
		h(A\cup B)\leq h(A)+h(B)\quad {\text{for all  }}\  A, B\in{\rm{Fin}}(S);
		$$
		\item[(2)] $h$ is right-subinvariant, i.e.,
		$$
		h(A s)\leq h(A)\quad {\text{for all  }}\  s \in S\  {\text {   and  }}\  A\in{\rm{Fin}}(S);
		$$
		\item[(3)] $h$ is bounded on singletons, i.e., there is $M>0$ such that
		$$
		h(\{s\})\leq M\quad {\text{for all  }}\  s \in S.
		$$
	\end{itemize}
	Then there exists a real number $\lambda\geq 0$, depending only on $h$, such that the sequence $\left(h(F_n)/|F_n|\right)_{n\in \mathbb{N}}$
	converges to $\lambda$ for every left F{\o}lner sequence $(F_n)_{n\in \mathbb{N}}$ of $S$.
\end{theorem}

\subsection{Actions of semigroup.}

Let $(X,d)$ be a compact metric space.  An {\it {action}} $S\curvearrowright X$ is a family of continuous maps $x\mapsto sx$ satisfying
\[
s(tx)=(st)x
\qquad(s,t\in S,\ x\in X).
\]
The action is called \emph{surjective} if each map $x\mapsto sx$ is onto.

A {\it dynamical system} is a pair $(X, S)$ consisting of a compact metrizable space $X$, called the {\it {phase space}}, equipped with a continuous action of a
countable semigroup $S$.

For a finite set $F\subset S$, define the Bowen pseudometric
\[
d_F(x,y)=\max_{s\in F}d(sx,sy),
\]
with the convention $d_\varnothing\equiv0$.

For any subsets $Z\subseteq X$ and $K\subseteq S$, and any $s\in S$, we write
\begin{equation*}
	s Z=\{s x  : x\in Z\},\quad K x=\{s x  : s\in K\},\quad K Z=\{s x : x\in Z, s\in K\}.
\end{equation*}
The $S$-{\it {orbit}} of a point $x\in X$ is the set $Sx$. A set $Z\subseteq X$ is said to be $S$-{\it {invariant}} if $S A\subseteq A$.


\subsection{Topological entropy of semigroup action}

\subsubsection{the definition of an open cover}
Given two open covers \(\mathscr{U}\) and \(\mathscr{V}\) of \(X\), their {\emph{join}} is defined by
\[
\mathscr{U}\vee\mathscr{V}=\{U\cap V : U\in\mathscr{U},\; V\in\mathscr{V}\},
\]
which is again an open cover of \(X\). Similarly one defines the join \(\mathscr{U}_1\vee\cdots\vee\mathscr{U}_n\) of finitely many open covers. For an element \(s\in S\) (which acts as a continuous map on \(X\)) and an open cover \(\mathscr{U}\), we set
\[
s^{-1}\mathscr{U}=\{s^{-1}(U) : U\in\mathscr{U}\},
\]
which is again an open cover. If \(\mathscr{U}\) is a finite open cover and \(F\subset S\) is a nonempty finite subset, we write
\[
\mathscr{U}^F = \bigvee_{s \in F} s^{-1}\mathscr{U}.
\]

If \(\mathscr{U}\) and \(\mathscr{V}\) are open covers of \(X\) and for every \(V\in\mathscr{V}\) there exists some \(U\in\mathscr{U}\) with \(V\subseteq U\), we say that \(\mathscr{V}\) is a \emph{refinement} of \(\mathscr{U}\) and denote this by \(\mathscr{V}\succeq\mathscr{U}\) or \(\mathscr{U}\preceq\mathscr{V}\). In particular, if every member of \(\mathscr{V}\) already belongs to \(\mathscr{U}\), then \(\mathscr{V}\) is a \emph{subcover} of \(\mathscr{U}\); obviously \(\mathscr{V}\succeq\mathscr{U}\) in this case.

\begin{definition}
	For an open cover \(\mathscr{U}\) of \(X\), let \(N(\mathscr{U})\) denote the minimal cardinality of a subcover of \(\mathscr{U}\), i.e.
	\[
	N(\mathscr{U}) = \min \Bigl\{ |\mathscr{V}| : \mathscr{V} \subset \mathscr{U},\; \bigcup_{V\in\mathscr{V}} V = X \Bigr\}.
	\]
\end{definition}

For a fixed open cover \(\mathscr{U}\) of \(X\) define a map \(h_{\mathscr{U}}\colon \rm{Fin}(S)\to \mathbb{R}\) by
\[
h_{\mathscr{U}}(F) = \log N\bigl(\mathscr{U}^F\bigr).
\]
By Proposition~5.2 of \cite{CCK}, the function \(h_{\mathscr{U}}\) is nonnegative, subadditive, right-invariant and bounded on singletons. Hence Theorem~\ref{fekete} immediately yields the following result, which is \cite[Theorem~5.3]{CCK}.

\begin{theorem}[{\cite{CCK}}]\label{thm:top-ent-cover}
	Let \(X\) be a compact metrizable space, \(S\) a countable cancellative left amenable semigroup, and \(\mathscr{U}\) an open cover of \(X\). Then for every left F{\o}lner sequence \(\{F_n\}_{n\in\mathbb{N}}\) in \(S\) the limit
	\[
	h_{\mathrm{top}}(\mathscr{U}) := \lim_{n\to\infty} \frac{h_{\mathscr{U}}(F_n)}{|F_n|}
	\]
	exists, is finite, and is independent of the chosen F{\o}lner sequence.
\end{theorem}

With this at hand, we define the topological entropy of the system.

\begin{definition}\label{def:top-ent}
	The \emph{topological entropy} of \((X, S)\) is
	\[
	h_{\mathrm{top}}(X) = \sup_{\mathscr{U}} h_{\mathrm{top}}(\mathscr{U}),
	\]
	where the supremum is taken over all finite open covers of \(X\).
\end{definition}

\subsubsection{Bowen's definition}

We now introduce Bowen's definition of topological entropy. For a large class of systems, Bowen's definition makes the computation of topological entropy more tractable. Our exposition in this part mainly follows \cite{KL}.

\begin{definition}
	Let \((X, S)\) be a dynamical system, \(d\) a metric on \(X\) compatible with the topology, \(\epsilon>0\), and \(Z\subseteq X\) a closed set.
	\begin{itemize}
		\item[(1)] A subset \(D\subset Z\) is called \((d, \epsilon)\)-{\it{separated}} for \(Z\) if \(d(x, y)\ge \epsilon\) for all distinct \(x, y\in D\).
		\item[(2)] A subset \(E\subset X\) is called \((d, \epsilon)\)-{\it{spanning}} for \(Z\) if for every \(x\in Z\) there exists \(y\in E\) with \(d(x, y)<\epsilon\).
	\end{itemize}
\end{definition}

\begin{definition}
	Let \(F\subset S\) be a nonempty finite set, \(\epsilon>0\), and \(Z\subseteq X\).
	\begin{itemize}
		\item[(1)] A subset \(D\subset Z\) is called a \((d_F, \epsilon)\)-{\it{separated set}} for \(Z\) if \(d_F(x, y)\ge \epsilon\) for all distinct \(x, y\in D\).
		\item[(2)] A subset \(E\subset X\) is called a \((d_F, \epsilon)\)-{\it{spanning set}} for \(Z\) if for every \(x\in Z\) there exists \(y\in E\) with \(d_F(x, y)<\epsilon\).
	\end{itemize}
\end{definition}

\begin{definition}\label{def:sep-span-card}
	For \(F\subset S\) nonempty finite, \(\epsilon>0\) and \(Z\subseteq X\), set
	\[
	\begin{aligned}
		\operatorname{sep}(d_F, \epsilon, Z)   & = \max \bigl\{ |D| : D \text{ is a } (d_F, \epsilon)\text{-separated set for } Z \bigr\}, \\
		\operatorname{span}(d_F, \epsilon, Z)  & = \min \bigl\{ |E| : E \text{ is a } (d_F, \epsilon)\text{-spanning set for } Z \bigr\}.
	\end{aligned}
	\]
\end{definition}

Now fix a left F{\o}lner sequence \(\mathcal{F}=\{F_n\}_{n\in\mathbb{N}}\) in \(S\). We define two metric entropies by
\begin{align*}
	h_{\mathrm{sep}}^{\mathcal{F}}(d, \epsilon, Z) & = \limsup_{n \to \infty} \frac{\log \operatorname{sep}(d_{F_n}, \epsilon, Z)}{|F_n|}, \\
	h_{\mathrm{span}}^{\mathcal{F}}(d, \epsilon, Z) & = \limsup_{n \to \infty} \frac{\log \operatorname{span}(d_{F_n}, \epsilon, Z)}{|F_n|},
\end{align*}
and then let
\begin{align*}
	h_{\mathrm{sep}}^{\mathcal{F}}(d, Z) & = \sup_{\epsilon > 0} h_{\mathrm{sep}}(d, \epsilon, Z), \\
	h_{\mathrm{span}}^{\mathcal{F}}(d, Z) & = \sup_{\epsilon > 0} h_{\mathrm{span}}(d, \epsilon, Z).
\end{align*}

\begin{example}[Dependence on the choice of F{\o}lner sequence]
	The above entropies of a non-invariant closed subset may depend on the choice of the F{\o}lner sequence, even for the full shift over \({\mathbb {Z}}\).
	
	Let \(S=\mathbb {Z}\), written additively, and consider the full shift
	\[
	X=\{0,1\}^{\mathbb{Z}}
	\]
	equipped with the shift action
	\[
	(T^{n} x)_k=x_{k+n},
	\qquad x\in X,\  n,k\in\mathbb{Z}.
	\]
	We endow \(X\) with the compatible metric
	\[
	d(x,y)=\sum_{k\in\mathbb {Z}}2^{-|k|}|x_k-y_k|.
	\]
	In particular, if \(x_0\neq y_0\), then \(d(x,y)\geq 1\).
	
	We shall construct a closed subset \(Z\subset X\) and two F{\o}lner sequences
	\(\mathcal{F}=\{F_n\}_{n\geq 1}\) and \(\mathcal{G}=\{G_n\}_{n\geq 1}\) such that
	\[
	h_{\mathrm{sep}}^{\mathcal{F}}(d, Z)=\log 2
	\quad\text{and}\quad
	h_{\mathrm{sep}}^{\mathcal{G}}(d, Z)=0.
	\]
	
	Set
	\[
	b_0=0,
	\qquad
	b_m=\sum_{j=1}^{m}j^2,
	\]
	and define
	\[
	I_m=\{b_{m-1}, b_{m-1}+1,\cdots, b_m-1\}.
	\]
	Then
	\[
	|I_m|=m^2.
	\]
	Let
	\[
	A=\bigcup_{m\geq 1;\, m\  {\text {is odd}}} I_m
	\]
	and define
	\[
	Z=Z_A:=
	\bigl\{
	x\in X \,:\, x_k=0\quad \text{for every } \quad k\notin A
	\bigr\}.
	\]
	Thus the coordinates indexed by \(A\) are free, whereas all coordinates outside
	\(A\) are fixed to be \(0\). Since
	\[
	Z_A=\bigcap_{k\notin A}\{x\in X \, :\, x_k=0\},
	\]
	the set \(Z_A\) is closed, and hence compact. Notice that \(Z_A\) is not \(T\)-invariant.
	
	We first record an estimate for the separated numbers of \(Z_A\). Let
	\(F\subset \mathbb {Z}\) be nonempty and finite. Then
	\begin{equation}\label{example01}
		\operatorname{sep}(d_F, 1, Z_A)\geq 2^{|A\cap F|}.
	\end{equation}
	Indeed, for each map
	\[
	u : A\cap F\longrightarrow \{0,1\},
	\]
	define \(x^{u}\in Z_A\) by
	\[
	x^{u}_k=
	\begin{cases}
		u(k), & k\in A\cap F,\\
		0, & k\notin A\cap F.
	\end{cases}
	\]
	If \(u\neq v\), then \(u(s)\neq v(s)\) for some \(s\in A\cap F\). Hence
	\[
	(T^{s} x^{u})_0=x^{u}_s\neq x^{v}_s=(T^{s} x^{v})_0,
	\]
	and therefore
	\[
	d_F(x^{u}, x^{v})\geq d(T^{s} x^{u}, T^{s} x^{v})\geq 1.
	\]
	Consequently, the collection of all \(x^{u}\) is \((d_F, 1)\)-separated, which
	proves (\ref{example01}).
	
	Conversely, fix \(\varepsilon>0\). Choose \(r=r(\varepsilon)\in\mathbb {N}\)
	such that
	\[
	\sum_{|k|>r}2^{-|k|}<\varepsilon,
	\]
	and write
	\[
	K_r=\{-r, -r+1, \cdots, r\}.
	\]
	We claim that
	\begin{equation}\label{example02}
		\operatorname{sep}(d_F, \varepsilon, Z_A)
		\leq 2^{|A\cap(F+K_r)|},
	\end{equation}
	where
	\[
	F+K_r=\{s+k : s\in F,  k\in K_r\}.
	\]
	
	To see this, suppose that \(x, y\in Z_A\) agree on
	\(A\cap(F+K_r)\). Since both \(x\) and \(y\) vanish outside \(A\), they in fact
	agree on the whole set \(F+K_r\). Hence, for every \(s\in F\) and every
	\(|k|\leq r\),
	\[
	(T^{s} x)_k=x_{s+k}=y_{s+k}=(T^{s} y)_k.
	\]
	It follows that
	\[
	d(T^{s} x, T^{s} y)
	\leq
	\sum_{|k|>r}2^{-|k|}
	<\varepsilon
	\]
	for every \(s\in F\), and therefore
	\[
	d_F(x, y)<\varepsilon.
	\]
	Thus, on any \((d_F,\varepsilon)\)-separated subset of \(Z_A\), the restriction
	map
	\[
	x\longmapsto x|_{A\cap(F+K_r)}
	\]
	is injective. This proves (\ref{example02}).
	
	Now let \(\mathcal{H}=\{H_n\}_{n\geq 1}\) be any sequence of finite integer intervals with
	\(|H_n|\to\infty\). Since
	\[
	|(H_n+K_r)\setminus H_n|\leq 2r,
	\]
	we have
	\[
	|A\cap(H_n+K_r)|
	\leq
	|A\cap H_n|+2r.
	\]
	It follows from (\ref{example02}) that, for every \(\varepsilon>0\),
	\begin{equation}\label{example03}
		\begin{aligned}
			h_{\mathrm {sep}}^{\mathcal{H}}(d, \varepsilon, Z_A)
			&=
			\limsup_{n\to\infty}
			\frac{1}{|H_n|}
			\log \operatorname{sep}(d_{H_n}, \varepsilon, Z_A)
			&\leq
			\left(
			\limsup_{n\to\infty}
			\frac{|A\cap H_n|}{|H_n|}
			\right)\log 2.
		\end{aligned}
	\end{equation}
	On the other hand, (\ref{example01}) yields
	\begin{equation}\label{example04}
		h_{\mathrm{sep}}^{\mathcal{H}}(d, 1, Z_A)
		\geq
		\left(
		\limsup_{n\to\infty}
		\frac{|A\cap H_n|}{|H_n|}
		\right)\log 2.
	\end{equation}
	Combining (\ref{example03}) and (\ref{example04}), we obtain
	\begin{equation}\label{example05}
		h_{\mathrm{sep}}^{\mathcal{H}}(d, Z_A)=
		\left(
		\limsup_{n\to\infty}
		\frac{|A\cap H_n|}{|H_n|}
		\right)\log 2.
	\end{equation}
	
	We now choose
	\[
	F_n=I_{2n-1}
	\quad\text{and}\quad
	G_n=I_{2n}.
	\]
	Both \(\mathcal{F}=\{F_n\}\) and \(\mathcal{G}=\{G_n\}\) are F{\o}lner sequences in \(\mathbb {Z}\), since
	they consist of intervals whose lengths tend to infinity. More explicitly, for
	every fixed \(t\in\mathbb {Z}\),
	\[
	|(t+F_n)\backslash F_n|\leq 2|t|
	\]
	and hence
	\[
	\frac{|(t+F_n)\backslash F_n|}{|F_n|}
	\longrightarrow 0.
	\]
	The same argument applies to \(\mathcal{G}=\{G_n\}\).
	
	By construction,
	\[
	F_n\subset A,
	\]
	and therefore
	\[
	\frac{|A\cap F_n|}{|F_n|}=1.
	\]
	Consequently, (\ref{example05}) gives
	\[
	h_{\mathrm{sep}}^{\mathcal{F}}(d, Z_A)=\log 2.
	\]
	
	On the other hand,
	\[
	G_n\cap A=\varnothing,
	\]
	so that
	\[
	\frac{|A\cap G_n|}{|G_n|}=0.
	\]
	Again by (\ref{example05}),
	\[
	h_{\mathrm{sep}}^{\mathcal{G}}(d, Z_A)=0.
	\]
	
	We have therefore constructed a compact metric \(\mathbb {Z}\)-system, a closed
	non-invariant subset \(Z_A\subset X\), and two F{\o}lner sequences such that
	\[
	h_{\mathrm{sep}}^{\mathcal{F}}(d, Z_A)
	\neq
	h_{\mathrm{sep}}^{\mathcal{G}}(d, Z_A).
	\]
	Hence, for a closed subset which is not assumed to be invariant, the quantity
	\(h_{\mathrm{sep}}^{\mathcal{F}}(d, Z)\) may depend on the choice of the F{\o}lner sequence.
\end{example}

But for $S$-invariant close subset $Z\subseteq X$, the topological entropy of $Z$ defined via separated and spanning sets does not depend on
the choice of the F{\o}lner sequence. To prove that these two metric entropies of $X$ as defined in Definition \ref{def:sep-span-card} coincide with the topological entropy defined in the previous subsection, we need the following auxiliary results. They are semigroup-action analogues of the corresponding statements for group actions found in \cite{KL}; the ideas of the proofs are essentially taken from Chapter~9 of that monograph.

\begin{lemma}\label{lemma:sep-product}
	For every \(\epsilon>0\) there exists a constant \(M>0\) such that for any nonempty finite subsets \(E, F\subset S\),
	\[
	\operatorname{sep}(d_{EF}, \epsilon, X) \le M^{|EF\setminus F|} \cdot \operatorname{sep}\!\Bigl(d_F, \frac{\epsilon}{2}, X\Bigr).
	\]
\end{lemma}

\begin{proof}
	Since \(X\) is compact, the open cover \(\{B_{d}(x, \epsilon/2) : x \in X\}\) admits a finite subcover; denote it by \(\mathscr{U}\).
	Let \(D\) be an \((d_{EF}, \epsilon)\)-separated subset of \(X\) of maximal cardinality, and choose inside \(D\) a \((d_F, \epsilon/2)\)-separated set \(D'\) of maximal cardinality. For each \(x\in D'\) set
	\[
	D'_x = \{ y\in D : d_F(x, y) < \epsilon/2 \}.
	\]
	If \(y\neq z\) belong to \(D'_x\), then there must exist some \(s\in EF\setminus F\) such that \(sy\) and \(sz\) lie in different members of \(\mathscr{U}\). Indeed, if for every \(s\in EF\) the points \(sy, sz\) belonged to the same member of \(\mathcal{U}\), the diameter of that member would give \(d(sy,sz)<\epsilon\) for all \(s\in EF\), whence \(d_{EF}(y,z)<\epsilon\), contradicting the fact that \(D\) is \((d_{EF}, \epsilon)\)-separated. Consequently,
	\[
	|D'_x| \le |\mathscr{U}|^{|EF\setminus F|}.
	\]
	Since \(D = \bigcup_{x\in D'} D'_x\), we obtain
	\[
	|D| \le |\mathscr{U}|^{|EF\setminus F|}\cdot |D'|
	\le |\mathscr{U}|^{|EF\setminus F|}\cdot \operatorname{sep}\!\Bigl(d_F, \frac{\epsilon}{2}, X\Bigr).
	\]
	Taking \(M = |\mathscr{U}|\) completes the proof.
\end{proof}

\begin{lemma}\label{lemma:sep-span-ineq}
	For any nonempty finite \(F\subset S\), \(\epsilon>0\) and \(Z\subseteq X\),
	\[
	\operatorname{span}(d_F, \epsilon, Z) \le \operatorname{sep}(d_F, \epsilon, Z)
	\le \operatorname{span}\!\Bigl(d_F, \frac{\epsilon}{2}, Z\Bigr).
	\]
\end{lemma}

\begin{proof}
	Every maximal \((d_F, \epsilon)\)-separated subset of \(Z\) is also \((d_F, \epsilon)\)-spanning for \(Z\), which gives the first inequality.
	For the second, let \(D\) be a \((d_F, \epsilon)\)-separated set for \(Z\) of maximal size and let \(E\) be a \((d_F, \epsilon/2)\)-spanning set for \(Z\) of minimal size. For each \(x\in D\) pick \(y_x\in E\) with \(d_F(x, y_x)<\epsilon/2\), and define \(\varphi: D\to E\) by \(\varphi(x)=y_x\). If \(\varphi(x_1)=\varphi(x_2)=y\) for distinct \(x_1,x_2\in D\), then
	\[
	d_F(x_1,x_2) \le d_F(x_1,y) + d_F(x_2,y) < \epsilon,
	\]
	which contradicts the separatedness of \(D\). Hence \(\varphi\) is injective, yielding \(|D|\le |E|\).
\end{proof}

\begin{lemma}\label{lemma:sep-span-vs-cover}
	Let \(d\) be a metric on \(X\), \(\mathscr{U}\) a finite open cover of \(X\), \(\epsilon>0\), and \(F\subset S\) a nonempty finite set.
	\begin{itemize}
		\item[(i)] If \(\epsilon\) is a Lebesgue number for \(\mathscr{U}\) with respect to \(d\), then
		\[
		N(\mathscr{U}^{F}) \le \operatorname{span}\!\Bigl(d_F, \frac{\epsilon}{2}, X\Bigr).
		\]
		\item[(ii)] If every member of \(\mathscr{U}\) has \(d\)-diameter less than \(\epsilon\), then
		\[
		\operatorname{sep}(d_F, \epsilon, X) \le N(\mathscr{U}^{F}).
		\]
	\end{itemize}
\end{lemma}

\begin{proof}
	(i) Let \(E\) be a \((d_F, \epsilon/2)\)-spanning set for \(X\) of minimal cardinality. Because \(\epsilon\) is a Lebesgue number, for each \(x\in E\) and \(s\in F\) we can find \(U_{x,s}\in\mathscr{U}\) with \(B_{d}(sx, \epsilon/2)\subseteq U_{x,s}\). For any \(y\in X\) choose \(x_y\in E\) such that \(d_F(x_y, y)<\epsilon/2\). Then \(s y\in B_{d}(s x_y, \epsilon/2)\subseteq U_{x_y, s}\) for all \(s\in F\), which implies
	\[
	y \in \bigcap_{s\in F} s^{-1}U_{x_y, s}.
	\]
	Thus \(\bigl\{\bigcap_{s\in F} s^{-1}U_{x,s} : x\in E\bigr\}\) is a subcover of \(\mathscr{U}^{F}\), and consequently
	\[
	N(\mathscr{U}^{F}) \le |E| = \operatorname{span}\!\Bigl(d_F, \frac{\epsilon}{2}, X\Bigr).
	\]
	
	(ii) Let \(\mathscr{V}\) be a subcover of \(\mathscr{U}^{F}\) of minimal cardinality and let \(D\) be a \((d_F, \epsilon)\)-separated set. For each \(x\in D\) choose \(V_x\in\mathcal{V}\) containing \(x\). If \(x_1\neq x_2\) in \(D\) and \(V_{x_1}=V_{x_2}=V\), then, because every member of \(\mathcal{U}\) has \(d\)-diameter less than \(\epsilon\), the set \(V\) has \(d_F\)-diameter less than \(\epsilon\), which forces \(d_F(x_1,x_2)<\epsilon\) and contradicts the choice of \(D\). Hence the assignment \(x\mapsto V_x\) is injective, and therefore \(|D|\le |\mathscr{V}| = N(\mathscr{U}^{F})\). Maximising over \(D\) yields the claimed inequality.
\end{proof}

\begin{theorem}\label{thm:entropy-equality}
	Let \((X, d)\) be a compact metric space and let \(S\) be a countable cancellative left amenable semigroup acting continuously on \(X\). Suppose that $\mathcal{F}=\{F_n\}_{n\geq 1}$ is a left F{\o}lner sequence of $S$. Then
	\[
	h_{\mathrm{top}}(X) = h^{\mathcal{F}}_{\mathrm{sep}}(d, X) = h^{\mathcal{F}}_{\mathrm{span}}(d, X).
	\]
\end{theorem}

\begin{proof}
	The equalities follow directly from Lemmas~\ref{lemma:sep-span-ineq} and~\ref{lemma:sep-span-vs-cover}.
\end{proof}

\begin{remark}
	The topological entropy of $(X, S)$ defined via separated and spanning sets does not depend on the choice of the F{\o}lner sequence.
\end{remark}

The next lemma is a standard fact. However, we have not been able to locate a precise statement and proof for semigroup dynamical systems in the literature; we therefore provide them here.

\begin{lemma}\label{lemma:sep-union}
	Let \(A_1,\dots,A_n\subset S\) be nonempty finite subsets, $s\in S$ and \(Z\subseteq X\).
	\begin{itemize}
		\item[(1)] For every \(\epsilon>0\),
		\[
		\operatorname{sep}\!\Bigl(d_{\cup_{i=1}^n A_i}, 2\epsilon, Z\Bigr)
		\le \prod_{i=1}^n \operatorname{sep}(d_{A_i}, \epsilon, Z).
		\]
		\item[(2)] For every \(\epsilon>0\), $K\in{\rm {Fin}}(S)$ and $s\in S$, one has
		\begin{equation*}
			\operatorname{sep}\!\Bigl(d_{K s}, \epsilon, Z\Bigr)
			\le \operatorname{sep}(d_{K}, \epsilon, s(Z)) \quad {\text {and}}\quad
			\operatorname{sep}\!\Bigl(d_{\{s\}}, \epsilon, Z\Bigr)
			\le \operatorname{sep}(d, \epsilon, s(Z)).
		\end{equation*}
		Moreover, if $s : X\rightarrow X$ is surjective, then
		\begin{equation*}
			\operatorname{sep}\!\Bigl(d_{K s}, \epsilon, X\Bigr)
			=\operatorname{sep}(d_{K}, \epsilon, X) \quad {\text {and}}\quad
			\operatorname{sep}\!\Bigl(d_{\{s\}}, \epsilon, X\Bigr)
			= \operatorname{sep}(d, \epsilon, X).
		\end{equation*}
	\end{itemize}
\end{lemma}

\begin{proof}
	\noindent (1) By Lemma~\ref{lemma:sep-span-ineq} it suffices to show
	\[
	\operatorname{sep}\!\Bigl(d_{\cup_{i=1}^n A_i}, 2\epsilon, Z\Bigr)
	\le \prod_{i=1}^n \operatorname{span}(d_{A_i}, \epsilon, Z).
	\]
	Let \(D\) be a \((d_{\cup_i A_i}, 2\epsilon)\)-separated subset of \(Z\) of maximal cardinality. For each \(i=1,\dots,n\) let \(E_i\) be a \((d_{A_i}, \epsilon)\)-spanning set for \(Z\) of minimal cardinality. For every \(x\in D\) pick \(y_i(x)\in E_i\) with \(d_{A_i}(x, y_i(x))<\epsilon\), and define \(\varphi: D\to \prod_{i=1}^n E_i\) by \(\varphi(x)=(y_1(x),\dots, y_n(x))\). If \(\varphi(x_1)=\varphi(x_2)\), then for each \(i\)
	\[
	d_{A_i}(x_1,x_2) \le d_{A_i}(x_1,y_i(x_1)) + d_{A_i}(x_2,y_i(x_2)) < 2\epsilon,
	\]
	which implies \(d_{\cup_i A_i}(x_1,x_2)<2\epsilon\), contradicting the separatedness of \(D\). Thus \(\varphi\) is injective, and
	\[
	|D| \le \prod_{i=1}^n |E_i|.
	\]
	This is precisely the required inequality.
	
	\noindent (2) It follows from the definition of separated set of dynamical system.
\end{proof}


\subsection{Weak specification property}

In this section, we introduce the weak specification property for dynamical systems. We also provide examples of systems that enjoy this property, as well as a simple example of a system that does not.

\begin{definition}\label{def:wsp}
	Let $(X,d)$ be a compact metric space.  An action $S\curvearrowright X$ is a family of continuous maps $x\mapsto sx$ satisfying
	\[
	s(tx)=(st)x
	\qquad(s,t\in S,\ x\in X).
	\]
	The action is called \emph{surjective} if each map $x\mapsto sx$ is onto.  A set $Y\subset X$ is $S$-invariant if $sY\subset Y$ for every $s\in S$.
	
	For a finite set $F\subset S$, define the Bowen pseudometric
	\[
	d_F(x,y)=\max_{s\in F}d(sx,sy),
	\]
	with the convention $d_\varnothing\equiv0$.
\end{definition}

Compactness allows one to use arbitrary orbit pieces, a fact needed in Section 5.

\begin{lemma}\label{lem:spec-infinite}
	Assume $X$ is compact and the action has weak specification in the sense of Definition~\ref{def:wsp}.  Fix $\varepsilon>0$ and a corresponding specification set $K$.  Let $(F_i)_{i\in I}$ be an arbitrary family of subsets of $S$ satisfying $KF_i\cap F_j=\varnothing$ for $i\neq j$, and let $(x_i)_{i\in I}\subset X$.  Then there exists $x\in X$ such that
	\[
	d(sx,sx_i)\le\varepsilon
	\qquad(i\in I,\ s\in F_i).
	\]
\end{lemma}

\begin{proof}
	Fix $\varepsilon>0$ and let $K\subset S$ be the corresponding finite set given by weak specification. Let $(F_i)_{i\in I}$ be an arbitrary family satisfying $KF_i\cap F_j=\varnothing$ for $i\ne j$, and let $(x_i)_{i\in I}\subset X$.
	
	Consider the collection of all constraints
	\[
	\mathcal T:=\{(i,s): i\in I,\ s\in F_i\}.
	\]
	For each finite subset $\mathcal C\subset \mathcal T$, define
	\[
	E_{\mathcal C}:=\{x\in X: d(sx,sx_i)\le \varepsilon,\ \forall (i,s)\in\mathcal C\}.
	\]
	
	We first show that $E_{\mathcal C}\ne\varnothing$. Fix a finite $\mathcal C$. Let $\{i_1,\dots,i_r\}$ be the distinct indices appearing in $\mathcal C$, and for each $j=1,\dots,r$, set
	\[
	G_j:=\{s\in F_{i_j}: (i_j,s)\in\mathcal C\}.
	\]
	We may ignore any empty $G_j$; assume all are nonempty. Since $G_j\subset F_{i_j}$, for $j\ne l$ we have
	\[
	K G_j\cap G_l \subset K F_{i_j}\cap F_{i_l}=\varnothing.
	\]
	Thus the finite family $G_1,\dots,G_r$ satisfies the separation condition required for weak specification. By weak specification, there exists $x_{\mathcal C}\in X$ such that for each $j$,
	\[
	d(sx_{\mathcal C},sx_{i_j})\le \varepsilon \qquad (s\in G_j).
	\]
	Hence $x_{\mathcal C}\in E_{\mathcal C}$, so $E_{\mathcal C}\ne\varnothing$.
	
	Next, the family $\{E_{\mathcal C}: \mathcal C\subset \mathcal T \text{ finite}\}$ has the finite intersection property: for any finite subcollection $E_{\mathcal C_1},\dots,E_{\mathcal C_n}$, let $\mathcal C=\bigcup_{k=1}^n \mathcal C_k$; then $\mathcal C$ is finite and $E_{\mathcal C}\subset E_{\mathcal C_k}$ for each $k$, so
	\[
	\bigcap_{k=1}^n E_{\mathcal C_k}\supset E_{\mathcal C}\ne\varnothing.
	\]
	
	Each $E_{\mathcal C}$ is closed: for fixed $s\in S$, the map $x\mapsto sx$ is continuous (since the action is continuous), hence $x\mapsto d(sx,sx_i)$ is continuous, so the set $\{x: d(sx,sx_i)\le\varepsilon\}$ is closed; $E_{\mathcal C}$ is a finite intersection of such sets.
	
	Since $X$ is compact, the family of nonempty closed sets with the finite intersection property has nonempty total intersection. Thus there exists
	\[
	x\in \bigcap_{\mathcal C\subset \mathcal T \text{ finite}} E_{\mathcal C}.
	\]
	Now take any $i\in I$ and $s\in F_i$. The singleton $\mathcal C=\{(i,s)\}$ is finite, so $x\in E_{\mathcal C}$, which gives $d(sx,sx_i)\le\varepsilon$. Since $i$ and $s$ were arbitrary, the desired point exists.
\end{proof}

We now present examples of systems satisfying the weak specification property.

\begin{example}
	Let \(\mathbb{N}_0 = \mathbb{N}\cup\{0\}\). Then the semigroup dynamical system \((\{0,1\}^{\mathbb{N}_0}, \mathbb{N}_0)\) has the weak specification property.
\end{example}

\begin{proof}
	Endow \(\{0,1\}^{\mathbb{N}_0}\) with the standard metric
	\[
	d(x,y) = \exp(-n), \qquad n = \inf\{m : x(m) \neq y(m)\}.
	\]
	From the definition of \(d\) it is immediate that for any \(\epsilon>0\) there exists \(N=N(\epsilon)\in\mathbb{N}\) such that, for all \(y,z\in\{0,1\}^{\mathbb{N}_0}\),
	\[
	y|_{[0,N]} = z|_{[0,N]} \quad\Longrightarrow\quad d(y,z) \le \epsilon.
	\]
	
	Set \(K = [0,N]\). Take an arbitrary family \(\{F_j\}_{j\in J}\) of subsets of \(\mathbb{N}_0\) satisfying
	\[
	(K+F_i) \cap F_j = \varnothing \qquad (i\neq j),
	\]
	and an arbitrary family \(\{x_j\}_{j\in J}\subset X\).
	
	We claim that
	\[
	(K+F_i) \cap (K+F_j) = \varnothing \qquad\text{for all } i\neq j.
	\]
	Indeed, if the intersection were nonempty, there would exist \(m,m'\in K\), \(g_i\in F_i\), \(g_j\in F_j\) with \(m+g_i = m'+g_j\). Assuming without loss of generality that \(m\le m'\), we would obtain \(g_i = (m'-m) + g_j\) with \(m'-m\in K\), contradicting the fact that \((K+F_i)\cap F_j = \varnothing\). The claim is thus proved.
	
	We now construct a point \(x\in\{0,1\}^{\mathbb{N}_0}\) as follows: for every \(j\in J\) and every \(m\in K+F_j\), set \(x(m)=x_j(m)\); for all remaining \(m\in\mathbb{N}_0\), set \(x(m)=0\).
	
	It remains to verify that \(x\) fulfills the condition required in the definition. For any \(j\in J\) and any \(m\in F_j\), the points \(\sigma^m x\) and \(\sigma^m x_j\) coincide on \(K\) because \(m+K\subseteq K+F_j\) and, by construction, \((\sigma^m x)(k)=x(m+k)=x_j(m+k)=(\sigma^m x_j)(k)\) for every \(k\in K\). Hence \(d(\sigma^m x, \sigma^m x_j)\le\epsilon\). As this holds for all \(m\in F_j\), we obtain
	\[
	d_{F_j}(x,x_j) = \max_{m\in F_j} d(\sigma^m x, \sigma^m x_j) \le \epsilon.
	\]
	Thus the weak specification property is satisfied.
\end{proof}

\begin{example}\label{wspNZ}
	Let $\mathbb N$ be equipped with the discrete metric
	\[
	\delta(a,b)=
	\begin{cases}
		0, & a=b,\\
		1, & a\neq b.
	\end{cases}
	\]
	Let
	\[
	X=\mathbb N^{\mathbb Z}
	\]
	and equip $X$ with the product metric
	\[
	d_X(x,y)=\sum_{k\in\mathbb Z} 2^{-|k|-1}\delta(x_k,y_k), \qquad x,y\in X.
	\]
	Consider the shift action of $\mathbb Z$ on $X$ given by
	\[
	(n\cdot x)_k=x_{k+n}, \qquad n,k\in\mathbb Z.
	\]
	Then the full shift system $(\mathbb N^{\mathbb Z},\mathbb Z)$
	satisfies the weak specification property.
\end{example}

\begin{proof}
	We use the following formulation of the weak specification property:
	
	for every $\varepsilon>0$, there exists a finite set $K\subset\mathbb Z$ such that, whenever
	$F_1,\ldots,F_r$ are finite subsets of $\mathbb Z$ satisfying
	\[
	(K+F_i)\cap F_j=\varnothing \qquad\text{for all }i\neq j,
	\]
	and $x^{(1)},\ldots,x^{(r)}\in X$ are arbitrary, there exists
	$y\in X$ such that
	\[
	d_X(s\cdot y,s\cdot x^{(i)})<\varepsilon
	\]
	for every $i\in\{1,\ldots,r\}$ and every $s\in F_i$.
	
	Let $\varepsilon>0$. Choose $L\in\mathbb N$ sufficiently large so that
	\[
	2^{-L}<\varepsilon.
	\]
	Set
	\[
	W=\{-L,-L+1,\ldots,L\}
	\]
	and define the finite buffer set
	\[
	K=W-W=\{-2L,-2L+1,\ldots,2L\}.
	\]
	
	Let $F_1,\ldots,F_r\subset\mathbb Z$ be finite sets satisfying
	\[
	(K+F_i)\cap F_j=\varnothing \qquad\text{whenever }i\neq j,
	\]
	and let
	\[
	x^{(1)},\ldots,x^{(r)}\in X
	\]
	be arbitrary prescribed points.
	
	We first show that the enlarged coordinate sets
	\[
	W+F_1,\ldots,W+F_r
	\]
	are pairwise disjoint. Suppose that, for some $i\neq j$, there exists
	\[
	q\in(W+F_i)\cap(W+F_j).
	\]
	Then there exist
	\[
	s\in F_i,\qquad t\in F_j,\qquad u,v\in W
	\]
	such that
	\[
	q=s+u=t+v.
	\]
	It follows that
	\[
	t=s+(u-v).
	\]
	Since $u-v\in W-W=K$, we have
	\[
	t\in K+F_i.
	\]
	Since also $t\in F_j$, this gives
	\[
	t\in(K+F_i)\cap F_j,
	\]
	contrary to the separation assumption. Therefore,
	\[
	(W+F_i)\cap(W+F_j)=\varnothing \qquad\text{for all }i\neq j.
	\]
	Fix an arbitrary symbol $a_*\in\mathbb N$. Define
	\[
	y=(y_q)_{q\in\mathbb Z}\in\mathbb N^{\mathbb Z}
	\]
	by
	\[
	y_q=
	\begin{cases}
		x^{(i)}_q,
		& \text{if }q\in W+F_i
		\text{ for some }i\in\{1,\ldots,r\},\\[1mm]
		a_*,
		& \text{if }q\notin
		\displaystyle\bigcup_{i=1}^{r}(W+F_i).
	\end{cases}
	\]
	Because the sets $W+F_i$ are pairwise disjoint, this definition is unambiguous.
	
	Fix $i\in\{1,\ldots,r\}$ and $s\in F_i$. If $|k|\leq L$, then $k\in W$, so
	\[
	s+k\in F_i+W.
	\]
	By the definition of $y$,
	\[
	y_{s+k}=x^{(i)}_{s+k}.
	\]
	Consequently,
	\[
	(s\cdot y)_k=y_{k+s}=x^{(i)}_{k+s}=(s\cdot x^{(i)})_k
	\]
	for every $k\in W$. Thus, the configurations $s\cdot y$ and
	$s\cdot x^{(i)}$ agree on all coordinates in $[-L,L]$.
	It follows that
	\begin{equation*}
		\begin{split}
			d_X(s\cdot y,s\cdot x^{(i)}) &=\sum_{k\in\mathbb Z} 2^{-|k|-1}\delta\bigl((s\cdot y)_k,(s\cdot x^{(i)})_k\bigr)\\
			&=\sum_{|k|>L} 2^{-|k|-1} \delta\bigl((s\cdot y)_k,(s\cdot x^{(i)})_k\bigr)\\
			&\leq\sum_{|k|>L}2^{-|k|-1}\\
			&=2^{-L}\\
			&<\varepsilon.
		\end{split}
	\end{equation*}
	Since $i$ and $s\in F_i$ were arbitrary, the point $y$ satisfies
	\[
	d_X(s\cdot y,s\cdot x^{(i)})<\varepsilon
	\]
	for every $i\in\{1,\ldots,r\}$ and every $s\in F_i$.
	
	Therefore, the full shift system
	\[
	(\mathbb N^{\mathbb Z},\mathbb Z)
	\]
	satisfies the weak specification property.
\end{proof}

The following generalization is obtained by an identical argument, with the fixed symbol $a_*$ replaced by any point of the alphabet.

\begin{example}\label{wspC}
	Let $(A, d_A)$ be a compact metric space, and let
	\[
	X=A^{\mathbb Z}.
	\]
	Consider the shift action of $\mathbb Z$ on $X$ defined by
	\[
	(n\cdot x)_k=x_{k+n}, \qquad n,k\in\mathbb Z,\quad x=(x_k)_{k\in\mathbb Z}\in X.
	\]
	Then the full shift system $(A^{\mathbb Z},\mathbb Z)$ has the weak
	specification property.
\end{example}

\begin{example}\label{nowsp}
	Let $S=\mathbb{Z}$ act trivially on $X=\{0, 1\}$; that is, for every $n\in \mathbb{Z}$, the map $n : X\to X$ is the identity map on $X$. The system $(X, S)$ does not have the weak specification property. Indeed, fix $\varepsilon<1/2$ and let $K\subset \mathbb{Z}$ be any finite set. Choose $n\in\mathbb{Z}$ sufficiently large so that $\{0\}$ and $\{n\}$ are $K$-separated in the sense that $K\cdot\{0\}\cap\{n\}=\varnothing$. Assign the target points $x_{\{0\}}=0$ and $x_{\{n\}}=1$. Since the action is trivial, any point $x\in X$ satisfies $sx=x$ for all $s\in S$. Thus $x$ cannot be within distance $\varepsilon<1/2$ of both $0$ and $1$ simultaneously. Hence the weak specification property fails.
\end{example}

\subsection{$m$-Expansivity and strong $h$-expansivity}

In this subsection we introduce \(m\)-expansive systems. Heuristically, an \(m\)-expansive system can be described as follows: there exists a constant such that any given orbit can stay close to at most \(m\) orbits (including itself) within that constant. Carvalho and Cordeiro \cite{CaCo} proved that for every \(m\in\mathbb{N}\) there exists a homeomorphism of a compact metric space possessing the shadowing property which is \(m\)-expansive but not \((m-1)\)-expansive. This shows that \(m\)-expansive systems form a class substantially larger than that of classical expansive systems, and therefore their investigation is of considerable interest.

Let $S$ be a countable semigroup, $(X, d)$ be a compact metric space and  $(X, S)$ be a dynamical system.

\begin{definition}
	The dynamical system $(X, S)$
	is said to be {\it {expansive}}
	if there exists $c>0$ with the property that if $x\not= y$ then
	$$
	\sup_{s\in S} d(s x, s y)>c.
	$$
	We call the constant $c$ an {\it {expansive constant}} for the $(X, S)$.
\end{definition}

Take $\delta>0$, we define
\begin{equation*}
	\Gamma_{\delta}(x, S):=\{y\in X\  |\  d(s x, s y)\leq\delta,\  \forall s\in S\}.
\end{equation*}
We often write $\Gamma_{\delta}(x)$ when the acting semigroup $S$ is clear from the context.

In \cite{Mora} Morales first introduced the notation of $m$-expansiveness, which is a natural generalization of the usual expansiveness.

\begin{definition}
	Let $m\in\mathbb{N}$, $S$ be a countable semigroup and $(X, d)$ be a compact metric space. A dynamical system $(X, S)$ is said to be
	$m$-{\it {expansive}} if there exists an $m$-{\it {expansive constant}} $c>0$ for $(X, S)$ such that for every $x\in X$, $\Gamma_{c}(x)$
	has at most $m$ elements.
\end{definition}

Since $\Gamma_{\delta}(x)$ is not $S$-invariant subset of $(X, S)$, we need to introduce a metric entropy of $\Gamma_{\delta}(x, S)$ is not
depend on the choice of F{\o}lner sequence of $S$.

\begin{definition}
	Let $S$ be a countable amenable semigroup and $(X, d)$ be a compact metric space. A dynamical system $(X, S)$ is said to be
	{\it {strong $h$-expansive}} if there exists an $h$-{\it {expansive constant}} $c>0$ for $(X, S)$ such that for every $x\in X$ and every
	left F{\o}lner sequence $\mathcal{F}=\{F_n\}_{n\geq 1}$ of $S$, one has
	$$
	h_{\rm {sep}}^{\mathcal{F}}(d, \Gamma_{c}(x))=0.
	$$
\end{definition}

It is clear that 1-expansiveness is just the classical expansiveness and $m$-expansiveness is strong $h$-expansiveness.

\section{A fixed-scale entropy theorem}
We now prove the technical result which replaces the standard fixed-scale entropy formula for expansive actions.

\begin{theorem}\label{expansive-entropy}
	Let $S$ be an infinite countable cancellative left-amenable semigroup acting continuously on a compact metric space $(X,d)$.  Assume that the action is strong $h$-expansive with constant $c>0$.  Let $\mathcal F=(F_n)$ be an increasing left F{\o}lner sequence with $\bigcup_nF_n=S$.  Then, for every $0<\varepsilon<c$,
	\[
	h_{\rm{top}}(X)
	=h_{\rm {sep}}^{\mathcal F}(d,\varepsilon,X)
	=h_{\rm {sep}}^{\mathcal F}(d,c,X).
	\]
	The same conclusion holds for every closed invariant subsystem $Y\subset X$, with the restricted action and the same F{\o}lner sequence.
\end{theorem}

\begin{proof}
	By Theorem \ref{foner_exist}, we know that such F{\o}lner sequence stated in theorem always exists. Fix $0<\varepsilon<c$.  Since separated numbers are monotone in the scale,
	\[
	h_{\rm {sep}}^{\mathcal F}(d,\varepsilon,X)
	\ge h_{\rm {sep}}^{\mathcal F}(d,c,X).
	\]
	It remains to prove the reverse inequality.

	Suppose that $O$ is a $(d_{F_n}, c)$-separated set for $X$ with maximum cardinality
	${\rm {sep}}(d_{F_n}, c, X)$ for $n\in\mathbb{N}$. Take $x\in O$ and denote by
	\begin{equation}\label{ee01}
		D_n(x, c):=\{y\in X : d_{F_n}(y, x)\leq c\}.
	\end{equation}
	
	Let $0<\varepsilon<\min(1/2, c)$ and let $k\in\mathbb{N}$. We take $\eta=\frac{\varepsilon}{2k}$. In what follows, we will estimate the upper bound of ${\rm {span}}(d_{F_n}, \varepsilon, D_n(x, c))$ for $n$ large enough.
	
	Since $S$ is a cancellative semigroup and  $0<\eta<1/2$, by Filling theorem (Theorem \ref{Filling theorem}), there exists an integer $N\geq 1$, depending only on $\eta$, such that the following holds.
	
	If $\{K'_i\}_{1\leq i\leq N}$ is a finite sequence of nonempty finite subsets of $S$ with
	\begin{equation}\label{ee02}
		\alpha(K'_j, K'_i)\leq \frac{\eta^{2 N}}{|K'_i|} \quad {\text {for all}}\quad 1\leq i<j\leq N
	\end{equation}
	and $F$ is a nonempty finite subset of $S$ such that
	\begin{equation}\label{ee03}
		\alpha(F, K'_i)\leq \eta^{2 N}  \quad {\text {for all}}\quad 1\leq i\leq N,
	\end{equation}
	then there exists a finite sequence $\{P_i\}_{1\leq i\leq N}$ of finite subsets of $S$ such that
	\begin{itemize}
		\item[(1)] each $P_i, 1\leq i\leq N$, is an $(\eta, K'_i)$-filling pattern of $F$;
		\item[(2)] the subsets $K'_iP_i\subset F$, $1\leq i\leq N$, are pairwise disjoint;
		\item[(3)] $F^*:=F\backslash \bigcup_{1\leq i\leq N} K'_iP_i$ satisfies $|F^*|\leq \eta |F|$.
	\end{itemize}
	
	Owing to $S$ being infinite, we can choose a finite subset $K\subset S$ with $|K|=2k$. Since $S$ is cancellative and left amenable, for above $\eta>0$, $N$ and $K\in{\rm {Fin}}(S)$, by (3) of Proposition \ref{infinite}, there exist $K_1, \cdots, K_N\in{\rm {Fin}}(S)$ satisfying
	\begin{itemize}
		\item[(1)] $K\subset K_i$ for all $1\leq i\leq N$;
		\item[(2)] $\alpha(K_i, K_j)\leq \frac{\eta^{2N}}{|K_j|}$, $1\leq j<i\leq N$.
	\end{itemize}
	
	Define
	\begin{equation}\label{K^*def}
		K^*:=\bigcup_{i=1}^{N}K_i.
	\end{equation}
	Since $(X, S)$ is a strong $h$-expansive system with expansive constant $c>0$, for any $y\in X$, one has
	\begin{equation*}
		h^{\mathcal{F}}_{\rm {sep}}(d, \Gamma_c(y))=0.
	\end{equation*}
	Owing to $\cup_{n\geq 1}F_n=S$ and the definition of $h^{\mathcal{F}}_{\rm {sep}}(d, \Gamma_c(y))$, for each $y\in X$, there exists
	finite subset $T(y)\in{\rm {Fin}}(S)$ satisfying that
	\begin{equation}\label{F_ydef}
		K^*\subset T(y)\quad {\text {and}}\quad {\rm {sep}}(d_{T(y)}, \varepsilon/12, \Gamma_c(y))<\exp{\frac{|T(y)|}{k}}.
	\end{equation}
	Let $E(y)$ be a $(d_{F(y)}, \varepsilon/12)$-separated set of $\Gamma_c(y)$ with cardinality ${\rm {sep}}(d_{F(y)}, \varepsilon/12, \Gamma_c(y))$,
	that is,
	\begin{equation}\label{E_ydef}
		|E(y)|<\exp{\frac{|F(y)|}{k}}\quad {\text {and}}\quad \Gamma_c(y)\subseteq \bigcup_{\omega\in E(y)} B_{d_{F(y)}}(\omega, \varepsilon/12).
	\end{equation}
	
	For each $y\in X$, we define
	\begin{equation}\label{Uydef}
		U(y):=\left\{z\in X : \exists\  \omega\in E(y), \  d_{F(y)}(z, \omega)<\frac{\varepsilon}{12}\right\}.
	\end{equation}
	By (\ref{E_ydef}), it is clear that
	\begin{equation*}
		\Gamma_c(y)\subseteq U(y).
	\end{equation*}
	Since $S$ is countable, we can represent $S$ as the union of an increasing finite subsets $S_l$ of $S$, that is, $S=\bigcup_{l=1}^{\infty}S_l$, $S_l\in{\rm {Fin}}(S)$ and $S_l\subset S_{l+1}$. Meanwhile, we know that
	\begin{equation}\label{contain}
		\bigcap_{l=1}^{\infty} \overline{B}_{S_l}(y, c)=\Gamma_c(y)\subseteq U(y).
	\end{equation}
	Thus there exists finite set $L(y)\in{\rm {Fin}}(S)$ such that
	\begin{equation*}
		\overline{B}_{d_{L(y)}}(y, c)\subseteq U(y).
	\end{equation*}
	Due to $U(y)$ being open, there is $\gamma_y>c$ such that
	\begin{equation}\label{Lydef}
		\overline{B}_{d_{L(y)}}(y, \gamma_y)\subseteq U(y).
	\end{equation}
	Let $V(y)$ be a open neighborhood of $y$ as follows:
	\begin{equation}\label{Vydef}
		V(y):=\{z\in X : d_{L(y)}(z, y)<\gamma_y-c\}.
	\end{equation}
	Thus $\{V(y) : y\in X\}$ forms an open cover of $X$. By compactness of $X$, there exists a finite subcover $V(y_1), \cdots, V(y_p)$, $p\in \mathbb{N}$, namely,
	$X\subseteq \bigcup_{l=1}^{p}V(y_l)$. Set
	\begin{equation}\label{Ldef}
		L=L(y_1)\cup\cdots\cup L(y_p)\in{\rm {Fin}}(S).
	\end{equation}
	Take
	\begin{equation}\label{etadef}
		\theta=\frac{\eta^{2 N}}{5|K^*|}.
	\end{equation}
	Since $S$ is a countable cancellative and left amenable semigroup and $\{F_n\}_{n\in\mathbb{N}}$ is a left F{\o}ner sequence of $S$, according to (2) of Proposition
	\ref{infinite}, we have $\lim_{n}\alpha(F_n, L)=0$ and $\lim_{n}\alpha(F_n, K_i)=0$ for all $1\leq i\leq N$. Then there exists $M_1\in\mathbb{N}$ such that, for all $n\geq M_1$, one has
	\begin{equation}\label{Intinequa}
		|F_n\backslash {\rm {Int}}_L(F_n)|<\theta |F_n| \quad {\text {and}}\quad  \alpha(F_n, K_i)<\theta\quad {\text {for all}}\quad 1\leq i\leq N.
	\end{equation}
	Noting ${\rm {Int}}_L(F_n)\subset F_n$ and by Proposition
	\ref{racp}, (\ref{etadef}) and (\ref{K^*def}), we deduce that
	\begin{equation*}
		\alpha({\rm {Int}}_L(F_n), K_i)\leq 4|K_i|\theta=4|K_i|\frac{\eta^{2 N}}{5|K^*|}<\eta^{2 N}.
	\end{equation*}
	This inequality shows that, when $n\geq M_1$, ${\rm {Int}}_L(F_n)$ satisfies the inequality (\ref{ee03}) for the finite subsets $K_1, \cdots, K_N$ of $S$. Hence by Filling Theorem as stated in the third paragraph in the proof of this theorem,  there exists a finite sequence $\{P_{i}\}_{1\leq i\leq N}$ of finite subsets of $S$ such that
	\begin{itemize}
		\item[(1)] each $P_i, 1\leq i\leq N$, is an $(\eta, K_i)$-filling pattern of ${\rm {Int}}_L(F_n)$;
		\item[(2)] the subsets $K_iP_i\subset {\rm {Int}}_L(F_n)$, $1\leq i\leq N$, are pairwise disjoint;
		\item[(3)] $F_n^*:={\rm {Int}}_L(F_n)\backslash \bigcup_{1\leq i\leq N} K_iP_i$ satisfies $|F_n^*|\leq \eta |{\rm {Int}}_L(F_n)|$.
	\end{itemize}
	
	\noindent {\bf {Claim 1.}}  One has $\sum_{i=1}^{N}|P_i|\leq \frac{|F_n|}{k}$.
	
	Given $1\leq i\leq N$, by the definition of $(\eta, K_i)$-filling parter, there exist $\widehat{K_ig}\subset K_ig$ for each $g\in P_i$ satisfying
	$\{\widehat{K_ig} : g\in P_i\}$ are pairwise disjoint and
	$$
	|\widehat{K_ig}|>(1-\eta)|K_ig|\geq \frac{|K_i|}{2}\geq \frac{|K|}{2}=k.
	$$
	Combining with the subsets $K_iP_i\subset {\rm {Int}}_L(F_n)$, $1\leq i\leq N$, being pairwise disjoint, we get
	\begin{equation*}
		|F_n|\geq\sum_{i=1}^{N}\sum_{g\in P_i}|\widehat{K_ig}|\geq k\sum_{i=1}^{N}|P_i|.
	\end{equation*}
	Hence Claim 1 is obtained.
	
	Recall that $O$ is a $(d_{F_n}, c)$-separated set for $X$, $x\in O$ and
	\begin{equation*}
		D_n(x, c):=\{y\in X : d_{F_n}(y, x)\leq c\}.
	\end{equation*}
	For each $(g, i)\in P_i\times \{i\}$ for $1\leq i\leq N$, since $\{V(y_1), \cdots, V(y_p)\}$ is an open cover of $X$, there exists
	\begin{equation}\label{Vdef}
		y_{(g, i)}\in \{y_1, \cdots, y_p\}\quad {\text {such that}}\quad g x\in V(y_{(g, i)}).
	\end{equation}
	
	\noindent {\bf {Claim 2.}}  If $g\in {\rm {Int}}_{L}(F_n)$ and $g x\in V\left(y_{(g, i)}\right)$, then
	\begin{equation*}
		g\left(D_n(x, c)\right)\subseteq U\left(y_{(g, i)}\right).
	\end{equation*}
	
	For any $v\in g\left(D_n(x, c)\right)$, there is some $u\in D_n(x, c)$ such that $g u=v$. Meanwhile, one has $d_{F_n}(u, x)\leq c$. The assumption
	$g\in {\rm {Int}}_{L}(F_n)$ means that $L g\subset F_n$, namely, $d_{L}(v, gx)\leq c$. The condition $g x\in V\left(y_{(g, i)}\right)$ shows that
	\begin{equation*}
		d_{L(y_{(g, i)})}(g x, y_{(g, i)})<\gamma_{y_{(g, i)}}-c.
	\end{equation*}
	By the definition of $L$ in (\ref{Ldef}), it follows that $L(y_{(g, i)})\subset L$. Hence we deduce that
	\begin{equation*}
		d_{L(y_{(g, i)})}(v, y_{(g, i)})<\gamma_{y_{(g, i)}},
	\end{equation*}
	that is, $v\in U\left(y_{(g, i)}\right)$. Therefore, Claim 2 is obtained.
	
	Now we denote by
	\begin{equation}\label{Rndef}
		R_n=F_n^*\cup (F_n\backslash {\rm {Int}}_{L}(F_n))=F_n^*\cup \partial_L(F_n).
	\end{equation}
	Fix $D$ is a $(d, \varepsilon/4)$-spanning set of $X$ and define index set $I$ as follows:
	\begin{equation}\label{geshu}
		I:=\prod_{(g, i)\in \bigcup_{j=1}^{N} P_j\times\{j\}} E(y_{(g, i)})\times \prod_{s\in R_n} D,
	\end{equation}
	where $E(y_{(g, i)})$ is defined as in (\ref{E_ydef}).
	For any point
	$$
	\mathbf{\overline{x}}=\left(x_{(g, i)}\right)_{(g, i)\in \bigcup_{j=1}^{N} P_j\times\{j\}}\times \left(x_s\right)_{s\in R_n}\in I,
	$$
	we define the subset $U(\mathbf{\overline{x}})$ of $D_n(x, c)$ as follows:
	\begin{equation}\label{Iset}
		\begin{split}
			U(\mathbf{\overline{x}})=\Big\{y\in D_n(x, c)\, :\, & d_{T(y_{(g, i)})}(g y, x_{(g, i)})<\varepsilon/6, \quad \forall\  (g, i)\in \bigcup_{j=1}^{N} P_j\times \{j\}; \\
			& d(s y, x_s)<\varepsilon/6,\quad \forall\  s\in R_n \Big\},
		\end{split}
	\end{equation}
	where $T(y_{(g, i)})$ is defined as in (\ref{F_ydef}).
	
	From (\ref{K^*def}) and (\ref{F_ydef}), we know that $K_i\subset K^*\subset T(y_{(g, i)})$ for $1\leq i\leq N$. Thus we have
	$\cup_{i=1}^{N}K_i P_i \cup R_n=F_n$. Hence it follows that, for any $y_1, y_2\in U(\mathbf{\overline{x}})$, we have $d_{F_n}(y_1, y_2)<\varepsilon/3$. Hence the
	$d_{F_n}$-diameter of $U(\mathbf{\overline{x}})$ is less than $\varepsilon/3$, namely,
	\begin{equation*}
		{\rm {diameter}}_{d_{F_n}}(U(\mathbf{\overline{x}}))\leq \varepsilon/3.
	\end{equation*}
	
	\noindent {\bf {Claim 3.}}  We have
	\begin{equation*}
		D_n(x, c)\subseteq \bigcup_{\mathbf{\overline{x}}\in I}U(\mathbf{\overline{x}}).
	\end{equation*}
	
	For any $y\in D_n(x, c)$ and any $(g, i)\in P_i\times \{i\}$, we know that there exists $y_{(g, i)}\in\{y_1, \cdots, y_p\}$ such that $g x\in V(y_{(g, i)})$ as in (\ref{Vdef}). Hence, by Claim 2, it follows that $g y\in U(y_{(g, i)})$.
	According to the definition of $U(y_{(g, i)})$ as in (\ref{Uydef}), there exists a point $x'_{(g, i)}\in E(y_{(g, i)})$ such that $d_{T(y_{(g, i)})}(x'_{(g, i)}, y_{(g, i)})<\varepsilon/12$.
	
	Meanwhile, owing to $D$ being a $(d, \varepsilon/2)$-spanning set of $X$, for each $s\in R_n$, there is $x'_s\in D$
	such that $d(s y, x'_s)<\varepsilon/6$. Take
	\begin{equation*}
		\mathbf{\overline{x'}}=\left(x'_{(g, i)}\right)_{(g, i)\in \bigcup_{j=1}^{N} P_j\times\{j\}}\times \left(x'_s\right)_{s\in R_n}.
	\end{equation*}
	Then $y\in U(\mathbf{\overline{x'}})$. Hence Claim 3 is obtained.
	
	Recall that $O$ is a $(d_{F_n}, c)$-separated set for $X$ with maximum cardinality
	${\rm {sep}}(d_{F_n}, c, X)$, $x\in O$ and the definition of $D_{n}(x, c)$ as in \ref{ee01}. Combining with Claim 3, it follows that
	\begin{equation*}
		X\subseteq \bigcup_{x\in E} D_{n}(x, c)\subseteq \bigcup_{x\in O} \bigcup_{\mathbf{\overline{x}}\in I}U(\mathbf{\overline{x}}).
	\end{equation*}
	We denote by
	$$
	T=T(y_1)\cup\cdots\cup T(y_p),
	$$
	where $T(y_i)$ is defined as in (\ref{F_ydef}).
	According to  the
	$d_{F_n}$-diameter of $O(\mathbf{\overline{x}})$ being less than $\varepsilon/3$ and Lemma \ref{lemma:sep-span-ineq}, we deduce that
	\begin{equation}\label{ee51}
		\begin{split}
			{\rm {sep}}(d_{F_n}, \varepsilon, X) & \leq {\rm {span}}(d_{F_n}, \varepsilon/2, X) \leq {\rm {sep}}(d_{F_n}, c, X)\cdot |I|  \\
			&  \leq {\rm {sep}}(d_{F_n}, c, X)\cdot \left(\prod_{(g, i)\in \cup_{j=1}^{N}P_j\times {j}}|E(y_{(g, i)})|\right)\cdot |D|^{|R_n|}
			\qquad ({\text {since (\ref{E_ydef})}})\\
			& \leq {\rm {sep}}(d_{F_n}, c, X)\cdot e^{\frac{|T|}{k}\cdot \sum_{i=1}^{N}|P_i|}\cdot |D|^{|R_n|}.
		\end{split}
	\end{equation}
	From the definitions of $R_n$ as in (\ref{Rndef}) and $F^*_n$, we have that
	\begin{equation*}
		\limsup_{n\rightarrow\infty}\frac{|R_n|}{|F_n|}\leq \limsup_{n\rightarrow\infty}\frac{|F^*_n|}{|F_n|}+\limsup_{n\rightarrow\infty}\frac{|\partial_L(F_n)|}{|F_n|}\leq \eta<\frac{1}{2k}.
	\end{equation*}
	Combining with (\ref{ee51}) and Claim 1, one has
	\begin{equation*}
		\limsup_{n\rightarrow\infty} \frac{\log{\rm {sep}}(d_{F_n}, \varepsilon, X)}{|F_n|} \leq \limsup_{n\rightarrow\infty} \frac{\log{\rm {sep}}(d_{F_n}, c, X)}{|F_n|}
		+ \frac{|T|}{k^2}+\frac{|D|}{2k}.
	\end{equation*}
	Letting $k\rightarrow\infty$, it follows that
	\begin{equation*}
		h^{\mathcal{F}}_{\mathrm{sep}}(d, \varepsilon, X)\leq h^{\mathcal{F}}_{\mathrm{sep}}(d, c, X),
	\end{equation*}
	and hence
	\begin{equation*}
		h^{\mathcal{F}}_{\mathrm{sep}}(d, \varepsilon, X)= h^{\mathcal{F}}_{\mathrm{sep}}(d, c, X).
	\end{equation*}
	Finally, Theorem \ref{thm:entropy-equality} gives
	\[
	h_{\rm {top}}(X)=h^{\mathcal{F}}_{\mathrm{sep}}(d, \varepsilon, X)=h_{\rm {sep}}^{\mathcal F}(d,c,X).
	\]
	
	The proof for a closed invariant subsystem is identical, since the restriction remains strong $h$-expansive with the same constant.
\end{proof}


\section{Proximal and asymptotic relations}

\begin{definition}
	A pair $(x,y)\in X^2$ is \emph{asymptotic} if
	\[
	d(sx,sy)\longrightarrow0
	\qquad(S\ni s\to\infty),
	\]
	meaning that for every $\varepsilon>0$ there is a finite $F\subset S$ such that $d(sx,sy)<\varepsilon$ for all $s\in S\setminus F$.  The asymptotic relation is an equivalence relation; its class through $x$ is denoted $\rm {Asy}(x)$.
	
	The pair $(x,y)$ is \emph{proximal} if
	\[
	\inf_{s\in S}d(sx,sy)=0.
	\]
	The proximal set of $x$ is
	\[
	{\rm {Prox}}(x)=\{y\in X:(x,y)\text{ is proximal}\}.
	\]
\end{definition}

The proximal relation need not be transitive.  Nevertheless, weak specification makes each proximal set large when the acting object has an identity.

\begin{proposition}\label{prop:proximal-dense-Gdelta}
	Let $S$ be an infinite left-cancellative monoid acting continuously on a compact metric space $(X,d)$.  If the action has the weak specification property, then ${\rm {Prox}}(y)$ is a dense $G_\delta$ subset of $X$ for every $y\in X$.
\end{proposition}

\begin{proof}
	Fix $y\in X$ and define
	\[
	U_n(y)=\bigcup_{s\in S}\{x\in X:d(sx,sy)<1/n\}.
	\]
	Each $U_n(y)$ is open and
	\[
	{\rm {Prox}}(y)=\bigcap_{n\ge1}U_n(y).
	\]
	It remains to show that $U_n(y)$ is dense.  Let $V\subset X$ be nonempty and open, choose $x_0\in V$, and choose $r>0$ with $B(x_0,r)\subset V$.  Fix
	\[
	0<\varepsilon<\min\{r,1/n\},
	\]
	and let $K\subset S$ be a specification set for $\varepsilon$.  Write $e$ for the identity of $S$.
	
	For $k\in K$, the set $\{t:kt=e\}$ contains at most one point by left cancellativity.  Since $S$ is infinite, we may choose $t\in S$ such that
	\[
	t\notin K,
	\qquad
	e\notin Kt.
	\]
	Thus the sets $F_1=\{e\}$ and $F_2=\{t\}$ satisfy
	\[
	KF_1\cap F_2=KF_2\cap F_1=\varnothing.
	\]
	Weak specification gives $z\in X$ with
	\[
	d(z,x_0)\le\varepsilon,
	\qquad
	d(tz,ty)\le\varepsilon.
	\]
	Hence $z\in V\cap U_n(y)$.  Therefore every $U_n(y)$ is dense, and the Baire category theorem completes the proof.	
\end{proof}

\begin{corollary}\label{0202}
	Under the hypotheses of Proposition~\ref{prop:proximal-dense-Gdelta}, for every $x,y\in X$ there exists $z\in X$ such that
	\[
	x,y\in{\rm {Prox}}(z).
	\]
\end{corollary}

\begin{proof}
	Both ${\rm {Prox}}(x)$ and ${\rm {Prox}}(y)$ are dense $G_\delta$ subsets of the compact metric space $X$, so their intersection is nonempty.  If $z\in {\rm {Prox}}(x)\cap {\rm {Prox}}(y)$, symmetry of proximality gives $x,y\in {\rm {Prox}}(z)$.
\end{proof}

\begin{corollary}
	Under the hypotheses of Proposition~\ref{prop:proximal-dense-Gdelta}, a map $f\colon X\to X$ is injective on every proximal set if and only if it is globally injective.
\end{corollary}

\begin{proof}
	Only the forward implication needs proof.  If $f(x)=f(y)$, choose $z$ from Corollary~\ref{0202}.  Since $x,y\in{\rm {Prox}}(z)$ and $f|_{{\rm {Prox}}(z)}$ is injective, $x=y$.
\end{proof}

The following finite-multiplicity lemma is the form of proximal rigidity needed for general semigroups, where an identity is not assumed.

\begin{lemma}\label{lem:proximal}
	Let $(X, d)$ be a compact metric space, $S$ be an infinite left-cancellative semigroup, \((X, S)\) be \(m\)-expansive with expansiveness constant \(c>0\). Suppose that \(\{x_1,\dots,x_{m+1}\}\) are \(m+1\) distinct points in \(X\) for which there exists \(F\in {\rm {Fin}}(S)\) such that
	\[
	d(s x_i, s x_j)\le c \qquad\text{for all } i\neq j \text{ and all } s\in S\backslash F.
	\]
	Then there exist distinct indices \(i',j'\in\{1,\dots,m+1\}\) such that \((x_{i'},x_{j'})\) is a proper proximal pair.
\end{lemma}

\begin{proof}
	We argue by contradiction. Assume that no two distinct points among \(\{x_1,\dots,x_{m+1}\}\) form a proximal pair. Then, for every \(i\neq j\),
	\[
	\inf_{s\in S} d(s x_i, s x_j) > 0.
	\]
	Consequently, there exists \(\kappa>0\) such that
	\begin{equation}\label{eq:separated}
		d(s x_i, s x_j)>\kappa \qquad\text{for all distinct } i,j.
	\end{equation}
	Since \(X\) is compact and $S$ is infinite, we can extract a pairwise distinct subsequence \(\{s_k\}_{k\in\mathbb{N}}\subset S\)  such that for every \(i\in\{1,\dots,m+1\}\) the limit
	\[
	\lim_{k\to\infty} s_k x_i = y_i
	\]
	exists in \(X\). Using \eqref{eq:separated} we obtain
	\[
	d(s_k x_i, s_k x_j)>\kappa
	\]
	for all distinct \(i,j\), and letting \(k\to\infty\) yields \(d(y_i, y_j)\ge\kappa\). Thus the points \(\{y_1,\dots,y_{m+1}\}\) are pairwise distinct.
	
	By the hypothesis there exists \(F\in {\rm {Fin}}(S)\) such that
	\[
	d(s x_i, s x_j)\le c \qquad\text{for all } i\neq j \text{ and all } s\in S\backslash F.
	\]
	Now fix an arbitrary \(s\in S\). By Fact \ref{finite set lemma}, there exists a finite set $T\in{\rm {Fin}}(S)$ such that
	\begin{equation*}
		\{s\}(S\backslash T)\subseteq S\backslash F.
	\end{equation*}
	Since $T$ is finite, there is $N\in\mathbb{N}$ satisfying that $s_k\in S\backslash T$ whenever $k\geq N$. Thus, for $k\geq N$, we have
	$(s s_k)\in S\backslash F$, and hence
	\[
	d(s s_k x_i, s s _k x_j)\le c \qquad\text{for all } i\neq j.
	\]
	Taking the limit \(k\to\infty\) and using the continuity of \(s : X\rightarrow X\) gives
	\[
	d(s y_i, s y_j)\le c.
	\]
	Since \(s\) was arbitrary, we conclude that
	\[
	\{y_1,\dots,y_{m+1}\}\subseteq \Gamma_c(y_1),
	\]
	which contradicts the fact that \((X,S)\) is \(m\)-expansive. Hence there must exist a proper proximal pair among the points \(\{x_1,\dots,x_{m+1}\}\).
\end{proof}


\section{Entropy gap and the proximal Myhill theorem}

The key ingredient for the proof of Theorem \ref{main} is the following classical result, which is of independent interest. It indicates that weak specification gives the $m$-expansive system a 'entropy irreducibility': the topological entropy of the true subsystem strictly decreases.

\begin{theorem}[Entropy gap]\label{thm:entropy-drop}
	Let $S$ be an infinite countable cancellative left-amenable semigroup acting continuously and surjectively on a compact metric space $(X,d)$.  Suppose that the action is $m$-expansive and has the weak specification property.  If $Y\subsetneq X$ is a nonempty closed $S$-invariant subset, then
	\[
	h_{\rm {top}}(Y)<h_{\rm {top}}(X).
	\]
\end{theorem}

\begin{proof}
	Let $c>0$ be an $m$-expansive constant and choose $z\in X\setminus Y$.  Set
	\[
	\eta=\min\left\{\frac{{\rm {dist}}(z,Y)}8,\frac c8\right\}>0.
	\]
	Apply weak specification with tracing error $\eta/2$, and let $K\subset S$ be a corresponding finite specification set.
	
	Fix an increasing exhaustive a left F{\o}lner sequence \(\mathcal{F}=\{F_n\}_{n\in\mathbb{N}}\), that is,
	$$
	F_n\subseteq F_{n+1}\quad {\text {and}}\quad \bigcup_{n=1}^{\infty} F_n=S.
	$$
	For each \(n\), let \(K_n\subseteq F_n\) be a subset of maximal cardinality satisfying
	\[
	K s \cap \{t\} = \varnothing \quad\text{for all distinct } s,t\in K_n.
	\]
	
	\noindent\textbf{Claim 1.} We have the inclusion
	\[
	F_n \subseteq K K_n \cup \Bigl(\bigcup_{k\in K} L_k^{-1}(K_n)\Bigr),
	\]
	where \(L_k\colon S\to S\) denotes left multiplication by \(k\).
	
	Indeed, suppose some \(t\in F_n\) does not belong to \(K K_n\) nor to any \(L_k^{-1}(K_n)\). Then for every \(s\in K_n\) we have \(K s\cap\{t\}=\varnothing\) because \(t\notin K K_n\). Moreover, if \(K t\cap K_n\neq\varnothing\), then \(t\in \bigcup_{k\in K}L_k^{-1}(K_n)\), contradicting the choice of \(t\). Hence \(K_n\cup\{t\}\) would still be a pairwise \(K\)-disjoint subset of \(F_n\), contradicting the maximality of \(K_n\). This establishes Claim~1.
	
	Since \(S\) is cancellative, each left translation \(L_k\) is injective; consequently \(|L_k^{-1}(K_n)|\le |K_n|\). Using Claim~1 we immediately obtain the estimate
	\begin{equation}\label{eq:Kn-size}
		|K_n| \ge \frac{|F_n|}{2|K|}.
	\end{equation}
	
	Now fix an arbitrary subset \(A\subseteq K_n\) and define
	\[
	L^{A}_K := \bigcup_{k\in K} L_k^{-1}(A).
	\]
	Because \(L_k\) is injective, \(|L^{A}_K|\le |K||A|\). Observe that
	\[
	F_n \subseteq K A \cup \bigl(F_n\setminus (K A\cup L^{A}_K)\bigr)\cup L^{A}_K.
	\]
	Applying Lemma~\ref{lemma:sep-union} we deduce
	\begin{equation}\label{eq:sep-estimate}
		\begin{split}
			\operatorname{sep}(d_{F_n}, 4\eta, Y)
			&\le \prod_{a\in A}\operatorname{sep}(d_{K a}, 2\eta, Y)
			\cdot \operatorname{sep}\bigl(d_{F_n\setminus (K A\cup L^{A}_K)}, 2\eta, Y\bigr)  \cdot \prod_{s\in L^{A}_K}\operatorname{sep}(d_{\{s\}}, 2\eta, Y)\\[4pt]
			&\le \bigl(\operatorname{sep}(d_K, 2\eta, X)\bigr)^{|A|}
			\cdot \operatorname{sep}\bigl(d_{F_n\setminus (K A\cup L^{A}_K)}, 2\eta, Y\bigr)
			\cdot \bigl(\operatorname{sep}(d, 2\eta, X)\bigr)^{|K||A|}\\[4pt]
			&= M^{|A|}\cdot \operatorname{sep}\bigl(d_{F_n\setminus (K A\cup L^{A}_K)}, 2\eta, Y\bigr),
		\end{split}
	\end{equation}
	where \(M = \operatorname{sep}(d_K, 2\eta, X)\cdot \bigl(\operatorname{sep}(d, 2\eta, X)\bigr)^{|K|}\).
	
	Now consider the family consisting of the singletons \(\{a\}\) for \(a\in A\) together with the set \(F_n\setminus (K A\cup L^{A}_K)\). We check that this family satisfies the \(K\)-disjointness condition required for the specification property. For this it suffices to verify
	\begin{equation}\label{eq:disjoint}
		K\bigl(F_n\setminus (K A\cup L^{A}_K)\bigr) \cap \{a\} = \varnothing \qquad (\forall\, a\in A).
	\end{equation}
	If \eqref{eq:disjoint} failed, there would exist \(a'\in A\) and \(f\in F_n\setminus (K A\cup L^{A}_K)\) with \(a' = k f\) for some \(k\in K\), which would yield \(f\in L_k^{-1}(a')\subseteq L^{A}_K\), contradicting the choice of \(f\). Hence \eqref{eq:disjoint} holds.
	
	Because the action of \(S\) on \(X\) is surjective, the map \(x\mapsto a x\) is onto for each \(a\in A\). For the point \(z\in X\setminus Y\) chosen at the beginning of the proof, pick \(z_a\in X\) with \(a z_a = z\). Let \(D_A\) be a maximal \((d_{F_n\setminus (K A\cup L^{A}_K)}, 2\eta)\)-separated subset of \(Y\). For any \(\omega\in D_A\) we apply the weak specification property to the family of points \(\{z_a\}_{a\in A}\cup\{\omega\}\) and the family of sets \(\{\{a\}\}_{a\in A}\cup\{F_n\setminus (K A\cup L^{A}_K)\}\). We obtain a point \(\omega_A\in X\) such that
	\[
	d_{\{a\}}(\omega_A, z_a) \le \frac{\eta}{2}\;\;(\forall\, a\in A)
	\quad\text{and}\quad
	d_{F_n\setminus (K A\cup L^{A}_K)}(\omega_A, \omega) \le \frac{\eta}{2}.
	\]
	Since \(a z_a = z\), this is equivalent to
	\begin{equation}\label{eq:spec-point}
		d(a\omega_A, z) \le \frac{\eta}{2}\;\;(\forall\, a\in A)
		\quad\text{and}\quad
		d_{F_n\setminus (K A\cup L^{A}_K)}(\omega_A, \omega) \le \frac{\eta}{2}.
	\end{equation}
	
	\noindent\textbf{Claim 2.} The set \(\{\omega_A : \omega\in D_A,\; A\subseteq K_n\}\) is \((d_{F_n}, \eta)\)-separated in \(X\).
	
	Take two distinct points \(\omega,\omega'\in D_A\) for the same \(A\). Using the definition of \(D_A\) and \eqref{eq:spec-point},
	\[
	\begin{split}
		d_{F_n}(\omega_A, \omega'_A)
		&\ge d_{F_n\setminus (K A\cup L^{A}_K)}(\omega_A, \omega'_A)\\
		&\ge d_{F_n\setminus (K A\cup L^{A}_K)}(\omega, \omega')
		- d_{F_n\setminus (K A\cup L^{A}_K)}(\omega_A, \omega)
		- d_{F_n\setminus (K A\cup L^{A}_K)}(\omega'_A, \omega')\\
		&\ge 2\eta - \tfrac{\eta}{2} - \tfrac{\eta}{2} = \eta.
	\end{split}
	\]
	
	Now consider distinct subsets \(A,B\subseteq K_n\) and points \(\omega\in D_A\), \(\omega'\in D_B\). Without loss of generality, there exists \(b\in B\setminus A\). From the maximality property of \(K_n\) one readily checks that \(b\in F_n\setminus (K A\cup L^{A}_K)\). Using \eqref{eq:spec-point} we have
	\[
	d(b\omega'_B, z)\le \frac{\eta}{2}, \qquad
	d(b\omega, b\omega_A)\le \frac{\eta}{2}.
	\]
	Since \(Y\) is invariant and \(\omega\in Y\), we have \(b\omega\in Y\). Therefore
	\[
	d(b\omega_A, b\omega'_B)
	\ge \mathrm{dist}(z,Y) - d(b\omega'_B, z) - d(b\omega, b\omega_A)
	\ge 8\eta - \frac{\eta}{2} - \frac{\eta}{2} \ge 7\eta \ge \eta.
	\]
	Consequently \(d_{F_n}(\omega_A, \omega'_B)\ge\eta\), which completes the proof of Claim~2.
	
	Observe that the cardinality of \(\{\omega_A\}\) is exactly \(\sum_{A\subseteq K_n} |D_A| = \sum_{A\subseteq K_n} \operatorname{sep}(d_{F_n\setminus (K A\cup L^{A}_K)}, 2\eta, Y)\). By Claim~2 and inequality \eqref{eq:sep-estimate},
	\begin{equation}\label{eq:count}
		\begin{split}
			\operatorname{sep}(d_{F_n}, \eta, X)
			&\ge \sum_{A\subseteq K_n} \operatorname{sep}\bigl(d_{F_n\setminus (K A\cup L^{A}_K)}, 2\eta, Y\bigr)\\
			&\ge \sum_{A\subseteq K_n} \operatorname{sep}(d_{F_n}, 4\eta, Y)\, M^{-|A|}.
		\end{split}
	\end{equation}
	
	Using the binomial formula,
	\[
	\sum_{A\subseteq K_n} t^{|A|} = \sum_{r=0}^{|K_n|} \binom{|K_n|}{r} t^r = (1+t)^{|K_n|} \qquad (t\in\mathbb{R}).
	\]
	Applying this with \(t = M^{-1}\) in \eqref{eq:count} yields
	\[
	\operatorname{sep}(d_{F_n}, \eta, X)
	\ge \operatorname{sep}(d_{F_n}, 4\eta, Y)\,(1+M^{-1})^{|K_n|}.
	\]
	
	Recall that \(4\eta < c/2\). By Theorem~\ref{expansive-entropy} and the estimate \eqref{eq:Kn-size},
	\[
	\begin{split}
		h_{\mathrm{top}}(X)
		&= \limsup_{n\to\infty} \frac{1}{|F_n|} \log \operatorname{sep}(d_{F_n}, \eta, X)\\
		&\ge \limsup_{n\to\infty} \frac{1}{|F_n|} \log \operatorname{sep}(d_{F_n}, 4\eta, Y)
		+ \liminf_{n\to\infty} \frac{|K_n|}{|F_n|} \log(1+M^{-1})\\[4pt]
		&\ge h_{\mathrm{top}}(Y) + \frac{1}{2|K|}\log(1+M^{-1})\\[4pt]
		&> h_{\mathrm{top}}(Y).
	\end{split}
	\]
	Thus \(h_{\mathrm{top}}(Y) < h_{\mathrm{top}}(X)\), as required.
\end{proof}

A direct consequence of the entropy drop theorem is the following corollary.

\begin{corollary}
	Let \((X,d)\) be a compact metric space, \(S\) a countable cancellative left amenable semigroup, and \((X,S)\) a dynamical system with surjective action that is expansive and has the weak specification property. If \(h_{\mathrm{top}}(X)=0\), then \(X\) contains no proper nonempty \(S\)-invariant closed subsets.
\end{corollary}

\begin{theorem}[Non-injectivity on proximal set]\label{thm:non-injective}
	Let $S$ be an infinite countable cancellative left-amenable semigroup acting continuously and surjectively on a compact metric space $(X,d)$.  Suppose that the action is $m$-expansive with constant $c>0$ and has the weak specification property.  Let $f\colon X\to X$ be continuous and $S$-equivariant.  If
	\[
	h_{\rm {top}}(f(X))<h_{\rm {top}}(X),
	\]
	then there exist distinct proximal points $x_1,x_2\in X$ such that
	\[
	f(x_1)=f(x_2).
	\]
\end{theorem}

\begin{proof}
	Let \(c>0\) be an expansiveness constant for \((X,S)\) and let \(\{F_n\}_{n\in\mathbb{N}}\) be a left F{\o}lner sequence. Set \(Y = f(X)\). Since \(h_{\mathrm{top}}(Y) < h_{\mathrm{top}}(X)\), there exists \(\eta>0\) such that \(h_{\mathrm{top}}(Y)+\eta < h_{\mathrm{top}}(X)\). Choose \(0<\epsilon < c/8\). By standard facts, \(h_{\mathrm{top}}(Y) = h_{\mathrm{sep}}(Y,\epsilon)\) and \(h_{\mathrm{top}}(X) = h_{\mathrm{sep}}(X,3\epsilon)\). From the definitions of these metric entropies we can ensure that for all sufficiently large \(n\),
	\[
	e^{\eta |F_n|}\,\mathrm{sep}(d_{F_n}, \epsilon, Y) < \mathrm{sep}(d_{F_n}, 3\epsilon, X).
	\]
	Let \(D\) and \(E\) be maximal \((d_{F_n}, 3\epsilon)\)-separated and \((d_{F_n}, \epsilon)\)-separated subsets of \(X\) and \(Y\) respectively. Clearly
	\[
	f(D) \subseteq \bigcup_{y\in E} B_{d_{F_n}}(y,\epsilon),
	\]
	so by the pigeonhole principle there exist \(y\in E\) and \(\Omega_n\subseteq D\) such that
	\begin{itemize}
		\item[(1)] \(|\Omega_n| \ge e^{\eta |F_n|}\);
		\item[(2)] \(f(\Omega_n) \subseteq B_{d_{F_n}}(y,\epsilon)\).
	\end{itemize}
	
	Uniform continuity of \(f\) provides \(0<\delta<\min\{\epsilon,\,c/3\}\) such that \(d(u,v)\le\delta\) implies \(d(f(u),f(v))<\epsilon\). Using the weak specification property with this \(\delta\), we find a finite set \(K\subseteq S\) such that for any family \(\{F_i\}_{i\in I}\) of subsets of \(S\) and any points \(\{x_i\}_{i\in I}\) with \(K F_i\cap F_j = \varnothing\) (\(i\neq j\)), there exists \(y\in X\) satisfying
	\[
	d_{F_i}(x_i, y) \le \delta \qquad (\forall\, i\in I).
	\]
	
	Consider the two sets \(E_1 = F_n\) and \(E_2 = S\setminus (K F_n\cup W)\), where
	\[
	W = \bigcup_{k\in K} L_k^{-1}(F_n).
	\]
	A straightforward verification shows \(K E_1\cap E_2 = \varnothing\) and \(K E_2\cap E_1 = \varnothing\). Fix an arbitrary point \(z\in X\). For each \(x\in \Omega_n\) we apply the specification property to the families \(\{E_1, E_2\}\) and \(\{x, z\}\) to obtain a point \(x^*\in X\) such that
	\begin{equation}\label{eq:star}
		d_{E_1}(x^*, x) \le \delta \quad\text{and}\quad d_{E_2}(x^*, z) \le \delta < c/2.
	\end{equation}
	Set \(\Omega^*_n = \{x^* : x\in \Omega_n\}\). The uniform continuity of \(f\) together with the choice of \(\delta\) gives
	\[
	d_{E_1}(f(x), f(x^*)) < \epsilon \quad\text{and}\quad d_{E_2}(f(x^*), f(z)) < \epsilon.
	\]
	For distinct \(x^*, \omega^*\in \Omega^*_n\) with corresponding \(x,\omega\in\Omega_n\), using condition (2) and \(E_1 = F_n\) we obtain
	\[
	d_{F_n}(f(x), y)<\epsilon,\;\; d_{F_n}(f(\omega), y)<\epsilon,\;\;
	d_{F_n}(f(x), f(x^*))<\epsilon,\;\; d_{F_n}(f(\omega), f(\omega^*))<\epsilon,
	\]
	and hence
	\[
	d_{F_n}(f(x^*), f(\omega^*)) < 4\epsilon.
	\]
	Thus for any two distinct points in \(\Omega^*_n\),
	\begin{equation}\label{eq:both}
		d_{F_n}(f(x^*), f(\omega^*)) < 4\epsilon \quad\text{and}\quad
		d_{E_2}(f(x^*), f(\omega^*)) < 2\epsilon.
	\end{equation}
	Moreover, since \(\Omega_n\) is \((d_{F_n}, 3\epsilon)\)-separated,
	\[
	d_{F_n}(x^*, \omega^*) \ge d_{F_n}(x,\omega) - d_{F_n}(x,x^*) - d_{F_n}(\omega,\omega^*) \ge 3\epsilon - \delta - \delta > \epsilon,
	\]
	which shows \(|\Omega^*_n| = |\Omega_n| \ge e^{\eta|F_n|}\). So we have found a set \(\Omega^*_n\) with the properties:
	\begin{enumerate}
		\item[(i)] \(|\Omega^*_n| \ge e^{\eta|F_n|}\);
		\item[(ii)] for all \(x^*, \omega^*\in\Omega^*_n\),
		\[
		d_{F_n}(f(x^*), f(\omega^*)) < 4\epsilon < c \quad\text{and}\quad
		d_{E_2}(f(x^*), f(\omega^*)) < 2\epsilon < c.
		\]
	\end{enumerate}
	
	Set \(P_n = (K F_n\cup W)\setminus F_n\), so that \(S = F_n \cup E_2 \cup P_n\). We estimate the size of \(P_n\). Clearly
	\[
	|P_n| \le |K F_n\setminus F_n| + |W\setminus F_n|
	\le \sum_{k\in K} |k F_n\setminus F_n| + \sum_{k\in K} |L_k^{-1}(F_n)\setminus F_n|.
	\]
	Because \(S\) is cancellative, each left translation \(L_k\) is injective and \(L_k(L_k^{-1}(F_n)\setminus F_n) = F_n\setminus k F_n\). Hence
	\[
	|L_k^{-1}(F_n)\setminus F_n| = |F_n\setminus k F_n| = |k F_n\setminus F_n|.
	\]
	Therefore
	\[
	|P_n| \le 2\sum_{k\in K} |k F_n\setminus F_n|.
	\]
	Since \(\{F_n\}\) is a F{\o}lner sequence, we obtain
	\begin{equation}\label{eq:Pn-small}
		\lim_{n\to\infty} \frac{|P_n|}{|F_n|} = 0.
	\end{equation}
	
	Now let \(L\) be a maximal \((d_{P_n}, \epsilon)\)-separated subset of \(Y\); then
	\[
	Y = \bigcup_{u\in L} B_{d_{P_n}}(u,\epsilon).
	\]
	Using Lemma~\ref{lemma:sep-union} and Theorem~\ref{expansive-entropy}, and the fact that \(\epsilon/2 < c/2\), we bound the cardinality of \(L\) by
	\begin{align*}
		|L| &\le \prod_{s\in P_n} \operatorname{sep}\bigl(d_{\{s\}}, \tfrac{\epsilon}{2}, Y\bigr)\\
		&\le \bigl(\operatorname{sep}(d, \tfrac{\epsilon}{2}, Y)\bigr)^{|P_n|}
		= \bigl(\operatorname{sep}(d, \tfrac{c}{2}, Y)\bigr)^{|P_n|}
		= M^{|P_n|},
	\end{align*}
	with \(M = \operatorname{sep}(d, c/2, Y)\). From \eqref{eq:Pn-small} it follows that for all large \(n\),
	\begin{equation}\label{eq:L-small}
		|L| \le M^{|P_n|} < e^{\frac{\eta}{3}|F_n|}.
	\end{equation}
	
	Applying the pigeonhole principle again, we can extract \(Q^*_n\subseteq \Omega^*_n\) and a point \(w\in L\) such that
	\begin{itemize}
		\item \(|Q^*_n| \ge \dfrac{|\Omega^*_n|}{|L|} \ge e^{\frac{\eta}{2}|F_n|}\);
		\item \(f(Q^*_n) \subseteq B_{d_{P_n}}(w,\epsilon)\).
	\end{itemize}
	In particular,
	\[
	d_{P_n}(f(x^*), f(\omega^*)) < 2\epsilon < c \qquad (\forall\, x^*,\omega^*\in Q^*_n).
	\]
	Combining this with property (ii) above and the decomposition \(S = F_n\cup E_2\cup P_n\), we obtain
	\[
	d_S(f(x^*), f(\omega^*)) < c \qquad (\forall\, x^*,\omega^*\in Q^*_n).
	\]
	Take a point $x^*_0\in Q^*_n$, one has
	\begin{equation*}
		f(Q^*_n)\subseteq \Gamma_{c}(f(x^*_0)).
	\end{equation*}
	Since the system $(X, S)$ is $m$-expansive with expansiveness constant \(c\), we have $|\Gamma_{c}(f(x^*_0))|\leq m$. Applying the pigeonhole principle, there exists
	$R^*_n\subseteq Q^*_n$ satisfying
	\begin{itemize}
		\item \(|R^*_n| \ge e^{\frac{\eta}{2}|F_n|}/m\);
		\item \(f(R^*_n)\) is a single point.
	\end{itemize}
	According to (1) of Proposition \ref{infinite}, we have $\lim_{n\rightarrow\infty} |F_n|=\infty$. This means that $|R^*_n|\geq m+1$ for $n$ large enough. Meanwhile, by (\ref{eq:star}), we know that
	\begin{equation*}
		d(s x^*, s \omega^*)<c \quad {\text {for all} }\quad x^*, \omega^*\in R^*_n\quad {\text {and}}\quad s\in S\backslash (K F_n\cup W).
	\end{equation*}
	By Lemma \ref{lem:proximal}, there exist two distinct points $x^*_1, x^*_2\in R^*_n$ which are proximal and $f(x^*_1)=f(x^*_2)$. The theorem follows.
\end{proof}


\section{A one-dimensional asymptotic version}

In this section $T\colon X\to X$ is a continuous map.  We use the forward Bowen ball
\[
\Gamma_c^+(x)=\{y:d(T^nx,T^ny)\le c\text{ for all }n\ge0\}.
\]
The map $T$ is positively $m$-expansive if $|\Gamma_c^+(x)|\le m$ for every $x$ and some $c>0$.

A sequence $(x_n)_{n\ge0}$ is a $\delta$-pseudo-orbit if
\[
d(Tx_n,x_{n+1})<\delta
\qquad(n\ge0).
\]
The map $T$ has the \emph{shadowing property} if for every $\varepsilon>0$ there exists $\delta>0$ such that every $\delta$-pseudo-orbit is $\varepsilon$-shadowed by a true orbit.

For $x\in X$ and $r>0$, set
\[
W_r^s(x)=\{y:d(T^nx,T^ny)\le r\text{ for every }n\ge0\}.
\]
Two points are forward asymptotic if $d(T^nx,T^ny)\to0$ as $n\to\infty$.

\begin{proposition}\label{prop:finite-tail-asymptotic-multiplicity}
	Let $T\colon X\to X$ be positively $m$-expansive with constant $c>0$ and suppose that $T$ has the shadowing property.  Then, for every $x\in X$,
	\[
	W_{c/4}^s(x)
	\]
	contains at most $m$ distinct forward asymptotic classes.
\end{proposition}

\begin{proof}
	Let \(x\in X\), and suppose $\{y_1, y_2,\cdots ,y_{m+1}\}\subseteq \mathcal{W}^{s}_{c/4}(x)$ are $m+1$ pairwise distinct points contained in $\mathcal{W}^{s}_{c/4}(x)$. It suffice to prove that there exist $1\leq i\neq j\leq m+1$ such that $y_i$ and $y_j$ are  asymptotic pair, that is,
	\begin{equation*}
		\lim_{n\to\infty} d\left(T^n(y_i), T^n(y_j)\right)=0.
	\end{equation*}
	
	We proceed by contradiction. Assume that for any distinct \(i,j\in\{1,\cdots, m+1\}\), the points \(y_i\) and \(y_j\) do not form an asymptotic pair. This means
	that there exists $r>0$ satisfies the following property: for any pair $(i, j)$ with $1\leq i\neq j\leq m+1,$ there exist a sequence of $\{a^{(i, j)}_n\}_{n\in\mathbb{N}}\subseteq \mathbb{N}$ such that
	that
	\begin{equation}\label{asy01}
		a^{(i,j)}_n \uparrow \infty\quad {\text {and}}\quad  d\left(T^{a^{(i,j)_{n}}}(y_i), T^{a^{(i,j)_{n}}}(y_j)\right)>r.
	\end{equation}
	
	For all distinct points pair $(i, j)$, $1\leq i, j\leq m+1$ and every $k\in\mathbb{N}$, we define:
	\begin{equation}\label{d_ij-def}
		d_{i, j}(k):=d\left(T^k(y_i), T^k(y_j)\right).
	\end{equation}
	Let
	\begin{equation}\label{eta-def}
		\eta:=\min_{i\neq j} \liminf_{k\to\infty} d_{i j}(k).
	\end{equation}
	
	\noindent {\it {Case 1.}}\quad  $\eta>0$.
	
	Then there exists $K\in\mathbb{N}$ such that
	\begin{equation}\label{asy02}
		d_{i j}(k)>\eta/2, \  \forall \ k\geq K, \  \forall i\neq j.
	\end{equation}
	By the compactness of \(X\), we may take a subsequence \(\{a_n\}_{n\in\mathbb{N}}\subseteq \mathbb{N}\) such that
	\begin{equation*}
		T^{a_n}(y_i)\rightarrow u_i \quad \forall i\in\{1, \cdots, m+1\}.
	\end{equation*}
	Fix $k\in \mathbb{N}\cup \{0\}$. Since $y_i, y_1\in \mathcal{W}^{s}_{c/4}(x)$, it is easy to see that
	\begin{equation*}
		d\left(T^k(T^{a_n}(y_i)), T^k(T^{a_n}(y_1))\right)=d\left(T^{k+a_n}(y_i), T^{k+a_n}(y_1)\right)<c/2.
	\end{equation*}
	Let $n\to \infty$, we get
	\begin{equation*}
		d(T^k(u_i), T^k(u_1))\leq c/2.
	\end{equation*}
	This means
	\begin{equation*}
		\{u_1, \cdots, u_{m+1}\}\subseteq \Gamma_c(u_1).
	\end{equation*}
	Meanwhile, the inequality (\ref{asy02}) shows that, for any $i\neq h$ and every $a_n\geq K$, one has
	$$
	d(T^{a_n}(y_i), T^{a_m}(y_j))=d_{ij}(a_n)>\eta/2.
	$$
	Let $n\to \infty$ again, we have
	$$
	d(u_i, u_j)\geq \eta/2>0.
	$$
	Thus $|\Gamma_c(u_1)|\geq m+1$ which contradict that $T$ is $m$-expansive.

	\noindent {\it {Case 2.}}\quad $\eta=0.$
	
	Then there exist $(u, v)$ with $1\leq u\neq v\leq m+1$ such that $\liminf_{k\to\infty} d_{u v}(k)=0$. Therefore, there exists a sequence of $\{b_n\}_{n\in\mathbb{N}}\subseteq \mathbb{N}$
	such that
	$$
	\lim_{n\to\infty} d_{u v}(b_n)=0.
	$$	
	In particular, we may choose a subsequence, still denoted by \(\{b_n\}_{n\in\mathbb{N}}\), such that:
	$$
	d_{u v}(b_n)<\frac{1}{n} \quad {\text {for all }}\  \forall n\in\mathbb{N}.
	$$
	Meanwhile, by (\ref{asy01}), we know that there is a sequence $\{a^{(u, v)}_n\}_{n\in\mathbb{N}}$ satisfying:
	$$
	d\left(T^{a^{(u,v)}_n}(y_u), T^{a^{(u,v)}_{n}}(y_v)\right)>r.
	$$
	We can construct two strictly increasing sequences of positive integers $\{c_n\}_{n\in\mathbb{N}}$ and $\{e_n\}_{n\in\mathbb{N}}$ by induction satisfying:
	\begin{itemize}
		\item[(1)] $\{c_n\}_{n\in\mathbb{N}}$ is a subsequence of $\{a^{(u,v)}_n\}_{n\in\mathbb{N}}$ and $d\left(T^{c_{n}}(y_u), T^{c_{n}}(y_v)\right)>r$ for all $n\in\mathbb{N}$;
		\item[(2)] $\{e_n\}_{n\in\mathbb{N}}$ is a subsequence of $\{b_n\}_{n\in\mathbb{N}}$ and $d_{u v}(e_n)<\frac{1}{n}$ for all $n\in\mathbb{N}$;
		\item[(3)] Interleaving condition: $e_n<c_n<e_{n+1}$ for every $n\in\mathbb{N}$.
	\end{itemize}
	Hence we have
	\begin{equation}\label{asy03}
		e_n<c_n<e_{n+1},\quad d_{u v}(c_n)>r,\quad d_{u v}(e_n)<\frac{1}{n}\quad {\text {for every  }}\  n\in\mathbb{N}.
	\end{equation}
	Take
	\begin{equation*}
		\alpha:=\min\left\{\frac{c}{4}, \frac{r}{4}\right\}.
	\end{equation*}
	By the shadowing property, there exists $\beta>0$ such that every $\beta$-pseudo-orbit can be $\alpha$-shadowed. Choose $\omega\in\mathbb{N}$ large enough, such that $1/\omega<\beta$.
	Now for every $l=1,2,\cdots, m+1$, we define a sequence $\{x^l_{k}\}_{k\in\mathbb{N}}$ as follows:
	\begin{equation}\label{asy04}
		x_k^l =
		\begin{cases}
			T^k(y_u), & k \le e_{w+l}-1, \\
			T^k(y_v), & k \ge e_{w+l}.
		\end{cases}
	\end{equation}
	
	Check the pseudo-orbit condition:
	\begin{itemize}
		\item At \(k = b_{w+l}\), we have
		\[
		d(T(x_{e_{\omega+l}-1}^l), x_{e_{\omega+l}}^l) = d(T^{e_{\omega+l}}(y_u), T^{e_{\omega+l}}(y_v)) < \frac{1}{\omega+l}< \frac{1}{\omega} < \beta.
		\]
		\item For other \(k\), the orbit is a genuine orbit, so the error is \(0\).
	\end{itemize}
	Therefore, \((x_k^l)_{k \in \mathbb{N}}\) is a \(\beta\)-pseudo-orbit.
	
	By shadowing property, there exist points $z^l\in X$ such that
	\begin{equation}\label{asy05}
		d\left(T^k(z_l), x_k^l\right) < \alpha, \quad \forall k \in \mathbb{N}.
	\end{equation}
	
	\noindent {\bf {Claim 1.}}\   $\{z_1,\cdots, z_{m+1}\}$ are pairwise disjoint.
	
	Let $1\leq l<s\leq m+1$. We will show that $z_l\neq z_s$. By (\ref{asy03}), we know that
	\begin{equation*}
		e_{\omega+l}<c_{\omega+s}\leq e_{\omega+s}-1.
	\end{equation*}
	Take $k'=c_{\omega+s}$. By (\ref{asy03}) and $k'\in\{c_n\}_{n\in\mathbb{N}}$, we have
	\begin{equation*}
		x^l_{k'}=T^{k'}(y_u), \  x^s_{k'}=T^{k'}(y_v)\quad {\text {and}}\quad d\left(x^l_{k'}, x^s_{k'}\right)=d_{u v}(k')>r.
	\end{equation*}
	By shadowing property, we know that
	\begin{equation*}
		d\left(T^{k'}(z_l), x^l_{k'}\right)<\alpha \quad {\text {and}}\quad d\left(T^{k'}(z_s), x^s_{k'}\right)<\alpha.
	\end{equation*}
	Thus,
	\begin{equation*}
		\begin{split}
			d\left(T^{k'}(z_l), T^{k'}(z_s)\right) & \geq  d\left(x^l_{k'}, x^s_{k'}\right)-d\left(T^{k'}(z_l), x^l_{k'}\right)-d\left(x^s_{k'}, T^{k'}(z_s)\right)\\
			& >r-2\alpha>r/2>0,
		\end{split}
	\end{equation*}
	which implies that $z_l\neq z_s$. Hence Claim is obtained.
	
	\noindent {\bf {Claim 2.}}\   $\{z_1,\cdots, z_{m+1}\}\subseteq \Gamma_c(z_1)$. 	
	
	For any $1\leq i\leq m+1$ and any $k\in\mathbb{N}$, we need to show that $d(f^k(z_i), f^k(z_1))\leq c$.
	
	By (\ref{asy04}) and the assumption, we know that
	$$
	z^1_k, z^i_k\in \{T^k(y_u), T^k(y_v)\}\quad {\text {and}}\quad d(T^(y_u), T^(y_v))\leq \frac{c}{4}.
	$$
	Thus we get
	$$
	d\left(z^1_k, z^i_k\right)\leq \frac{c}{4}.
	$$
	Hence, by (\ref{asy05}), we deduce that
	\begin{equation*}
		\begin{split}
			d\left(T^{k}(z_1), T^k(z_i)\right) & \leq  d\left(T^{k}(z_1), z^1_k\right)+ d\left(z^1_k, z^i_k\right) +d\left(z^i_k, T^k(z_i)\right)\\
			& \leq 2 \alpha+c/4<c,
		\end{split}
	\end{equation*}
	which implies that $z_i\in \Gamma_c(z_1)$. Therefore, Claim 2 is obtained.
	
	Claim 1 and Claim 2 show that $|\Gamma_c(z_1)|\geq m+1$ which contradicts the assumption that $T$ is $m$-expansive and hence the proposition is proved.
\end{proof}

\begin{theorem}\label{thm:asymptotic-collision}
	Let $T\colon X\to X$ be surjective, positively $m$-expansive, weakly specifying, and shadowing.  Let $f\colon X\to X$ be continuous and commute with $T$.  If
	\[
	h_{\rm {top}}(f(X))<h_{\rm {top}}(X),
	\]
	then $f$ is not injective on some forward asymptotic class.
\end{theorem}

\begin{proof}
	According to the Proof of Theorem \ref{thm:non-injective} as in (\ref{eq:star}), for large enough $n\in\mathbb N$, there exist
	$K_n\in\mathbb{N}$ and $R^*_n\subset\Omega^*_n\subset X$ satisfies:
	\begin{itemize}
		\item[(1)] for any $x^*_1, x^*_2 \in\Omega^*_n$ and all $k\geq K_n$, we have $d(T^k(x^*_1), T^k(x^*_2))<c/4$;
		\item[(2)] $|R^*_n|\geq m+1$;
		\item[(3)] $f(R^*_n)$ is a single point.
	\end{itemize}
	Suppose that $x^*_1, x^*_2,\cdots x^*_{m+1}\in R^*_n$ are $m+1$ pairwise disjoint points of $R^*_n$. Owing to (1) we know that, all $k\geq K_n$,
	\begin{equation*}
		d(T^k(x^*_i), T^k(x^*_j))<c/4\quad {\text {for all }}\quad k\geq K_n.
	\end{equation*}
	
	We claim that there exist $1\leq i\neq j\quad m+1$ such that $x^*_i, x^*_j$ are asymptotic pair.
	
	Let $y^*_i=T^{K_n}(x^*_i)$ for $1\leq i\leq m+1$. If there exists $1\leq i\neq j\leq m+1$ such that $y^*_i=y^*_j$, then $x^*_i, x^*_j$ are asymptotic pair,
	and hence the claim is  obtained. Thus we may assume that $\{y^*_1, y^*_2,\cdots y^*_{m+1}\}$ are pairwise disjoint. The conclusion
	(1) shows that
	\begin{equation*}
		\{y^*_1, y^*_2,\cdots y^*_{m+1}\}\subset \mathcal{W}^{s}_{c/4}(y^*_1),
	\end{equation*}
	where
	\[
	\mathcal{W}^{s}_{c/4}(y^*_1)
	:=
	\left\{
	y\in X:
	d(T^k (y^*_1),T^k (y))\leq c/4
	\text{ for every }k\geq 0
	\right\}
	\]
	By Proposition \ref{prop:finite-tail-asymptotic-multiplicity}, there exist $1\leq i\neq j\quad m+1$ such that $y^*_i, y^*_j$ are asymptotic pair, and hence
	$x^*_i, x^*_j$ are asymptotic pair. Then the restriction of \(f\) to some proximal class is not injective. Therefore, the proof of theorem is completed.
\end{proof}

\begin{proof}[Proof of Theorem~\ref{thm:one-dimensional-main}]
	If $f(X)\neq X$, Theorem~\ref{thm:entropy-drop}, applied to the $\mathbb N_0$-action generated by the surjective map $T$, gives
	\[
	h_{\rm {top}}(f(X))<h_{\rm {top}}(X).
	\]
	Theorem~\ref{thm:asymptotic-collision} then contradicts injectivity of $f$ on forward asymptotic classes.  Therefore $f$ is surjective.
\end{proof}


\section{Examples and sharpness phenomena}

The following examples isolate several hypotheses in Theorem~\ref{main}.

\begin{example}[Compactness cannot be dropped]\label{ex:noncompact}
	Let $S=\mathbb Z$ act by the shift on $X=\mathbb N^{\mathbb Z}$, where $\mathbb N$ has the discrete topology.  The shift is expansive, the action maps are surjective, and the full shift has weak specification.  Define
	\[
	f(x)_n=x_n+1.
	\]
	Then $f$ is continuous, equivariant, and globally injective, hence injective on every proximal set, but it is not surjective.  The missing hypothesis is compactness of $X$.
\end{example}

\begin{example}[$m$-expansiveness cannot be dropped]\label{ex:not-mexp}
	Let $S=\mathbb{Z}$, $X=[0, 1]^{\mathbb{Z}}$ and $(X, \mathbb{Z})$ be the Bernoulli shift. It is clear
	that $X$ is compact and for every $n\in \mathbb{Z}$, the map $n : X\to X$ is surjective. By Example \ref{wspC}
	we know that the system $(X, \mathbb{Z})$ has weak specification property. But for any $m\in\mathbb{N}$, the system $(X, \mathbb{Z})$ is not
	$m$-expansive. Define
	\begin{equation*}
		f : [0, 1]^{\mathbb{Z}}\to  [0, 1]^{\mathbb{Z}},\quad f(x)_n=\frac{x_n}{2}.
	\end{equation*}
	Then $f$ is continuous, $\mathbb{Z}$-equivariant and injective. But $f$ is not surjective because $f(X)\subset[0, 1/2]^{\mathbb{Z}}$.  Thus finite orbit multiplicity is essential for the theorem.
\end{example}

\begin{example}[Weak specification cannot be dropped]\label{ex:no-spec}
	Let $S=\mathbb Z$ act trivially on $X=\{0,1\}$ with the discrete metric.  The action is expansive and every action map is surjective, but the system does not have weak specification (see Example \ref{nowsp}).  The constant map
	\[
	f(0)=f(1)=0
	\]
	is continuous and equivariant.  Since the system is distal,
	\[
	{\rm {Prox}}(0)=\{0\},\qquad {\rm {Prox}}(1)=\{1\},
	\]
	so $f$ is injective on every proximal set, while $f$ is not surjective.
\end{example}

\begin{example}[Equivariance cannot be dropped]\label{ex:no-equivariance}
	Consider the Bernoulli shift $\mathbb Z\curvearrowright\{0,1\}^{\mathbb Z}$.  Choose a bijection
	\[
	b\colon\mathbb Z\longrightarrow\mathbb Z\setminus\{0\}
	\]
	and define
	\[
	f(x)_0=0,
	\qquad
	f(x)_{b(n)}=x_n.
	\]
	Then $f$ is continuous and injective, but not surjective.  It is not shift-equivariant.
\end{example}

\begin{example}[Infinitude of the acting group is essential]\label{ex:finite-group}
	Let
	\[
	G=C_2=\{e,a\},\qquad a^2=e,
	\]
	and let $X=\{0,1\}$ with $d(0,1)=1$.  Let $G$ act trivially on $X$ and define
	\[
	\tau(0)=0,\qquad \tau(1)=0.
	\]
	The action is expansive: $1/2$ is an expansive constant.  It also has the weak specification property in a degenerate sense.  Indeed, for every $\varepsilon>0$ choose $K=G$.  If a family $(F_i)$ satisfies
	\[
	KF_i\cap F_j=\varnothing\qquad(i\neq j),
	\]
	then at most one $F_i$ is nonempty, so the tracing condition is automatic.
	
	Since the action is trivial,
	\[
	\inf_{g\in G}d(g0,g1)=1,
	\]
	and hence
	\[
	{\rm {Prox}}(0)=\{0\},\qquad {\rm {Prox}}(1)=\{1\}.
	\]
	Thus $\tau$ is injective on every proximal set.  It is continuous and $G$-equivariant, but
	\[
	\tau(X)=\{0\}\neq X.
	\]
	This shows that the infinitude assumption in Theorem~\ref{main} cannot be omitted.
\end{example}

\begin{remark}
	The preceding examples do not address the optimality of amenability in the proximal formulation.  Unlike the classical homoclinic Garden of Eden theorem, proximal-set injectivity is closely related to global injectivity in weakly specifying systems, and the corresponding nonamenable problem is therefore entangled with surjunctivity phenomena.  We do not make a sharpness claim in that direction here.
\end{remark}


\end{document}